\documentclass[12pt]{amsart}
\usepackage{amsmath,amssymb,amsfonts,latexsym}
\usepackage{mathrsfs}
\usepackage[dvips]{geometry}
\usepackage{array}
\usepackage{epsf}
\usepackage{epsfig}
\newtheorem*{lemma}{Lemma}
\newtheorem*{prop}{Proposition}
\newtheorem*{thm}{Theorem}
\newtheorem*{cor}{Corollary}

\newcommand{\twoheaddownarrow}{\overset{\sim}{\twoheaddownarrow}}

\newcommand{\nc}{\newcommand}
\nc{\Ker}{\operatorname{Ker}} \nc{\rker}{\operatorname{rKer}}
\nc{\im}{\operatorname{Im}}
\nc{\stab}{\operatorname {Stab}}
\nc{\ann}{\operatorname {Ann}}
\nc{\Id}{\operatorname {Id}}
\nc{\Prim}{\operatorname {Prim}}
\nc{\Real}{\operatorname {Re}}
\nc{\Ext}{\operatorname {Ext}}
\nc{\rad}{\operatorname {rad}}
\nc{\rk}{\operatorname {rank}}
\nc{\Aut}{\operatorname {Aut}}
\nc{\supp}{\operatorname {supp}}
\nc{\height}{\operatorname {ht}}

\usepackage[utf8]{inputenc}
\usepackage[T1]{fontenc}
\usepackage{amsmath,amssymb,mathrsfs}
\usepackage{xcolor}
\usepackage{geometry}
\usepackage{caption}
\usepackage{tikz}
\usetikzlibrary{matrix}
\usepackage{graphicx}
\usepackage{float}
\usepackage{caption}

\usepackage{array}
\usepackage{epsf}
\usepackage{epsfig}
\usepackage{xcolor}

\usepackage{tikz,tikz-cd}

\usetikzlibrary{matrix}

\usetikzlibrary{fit}

\usetikzlibrary{matrix,backgrounds}
\pgfdeclarelayer{myback}
\pgfsetlayers{myback,background,main}

\tikzset{mycolor/.style = {line width=1bp,color=#1}}%
\tikzset{myfillcolor/.style = {draw,fill=#1}}%
\usetikzlibrary{arrows,matrix,positioning}
\usepackage{pstricks}

\usetikzlibrary{graphs,graphs.standard}
\usepackage{changepage}
\usepackage{tikz-cd}
\usepackage{mathrsfs}

\makeatletter
\newcommand*{\encircled}[1]{\relax\ifmmode\mathpalette\@encircled@math{#1}\else\@encircled{#1}\fi}
\newcommand*{\@encircled@math}[2]{\@encircled{$\m@th#1#2$}}
\newcommand*{\@encircled}[1]{%
  \tikz[baseline,anchor=base]{\node[draw,circle,outer sep=0pt,inner sep=.2ex] {#1};}}
\usepackage{float}

\begin{document}

\title [Reverse Tableaux]{The Reverse Tableaux: a Gateway to the Surjectivity of the Component Map.}
\author [Yasmine Fittouhi and Anthony Joseph]{Yasmine Fittouhi and Anthony Joseph}
\date{\today}
\maketitle
\vspace{-.9cm}\begin{center}

Department of Mathematics\\
The Weizmann Institute of Science\\
Rehovot, 7610001, Israel\\
fittouhiyasmine@gmail.com
\end{center}\

\

\

\vspace{-.9cm}\begin{center}
Donald Frey Professional Chair\\
Department of Mathematics\\
The Weizmann Institute of Science\\
Rehovot, 7610001, Israel\\
anthony.joseph@weizmann.ac.il
\end{center}\

\

\date{\today}
\maketitle

Key Words: Invariants, Parabolic adjoint action, Component map, Reverse Tableaux.

AMS Classification: 17B35

 \

\textbf{Abstract}.

Let $G$ be a simple algebraic group over the complex field $\mathbb C$, $B$ a fixed Borel subgroup, $P$ a parabolic subgroup containing $B$, $P'$ its derived group acting on the Lie algebra  $\mathfrak m$ of its nilradical.

The nilfibre $\mathscr N$ (for this action) is the zero locus  of the augmentation $\mathcal I_+$ of the semi-invariant algebra $\mathcal I:=\mathbb C[\mathfrak m]^{P'}$.   Via Richardson's theorem, $\mathcal I$ is polynomial.

Assume $G$ of type $A$.

Then the generators of $\mathcal I$ may be taken to be the Benlolo-Sanderson invariants \cite {BS}.

Following this, in Y. Fittouhi and A. Joseph, The Magic and Mystery of Component Tableaux, Indag Math., online January 2026, a set $\{\mathscr T^\mathcal C\}$ of component tableaux was constructed each encoding explicit combinatorial data $\mathcal C$.  Each tableau $\{\mathscr T^\mathcal C\}$ defines a component $\mathscr C$ of $\mathscr N$ and the resulting map $\{\mathscr T^\mathcal C\}\mapsto \mathscr C$  was shown to be injective.

Here this data is more simply encoded in a multi-set, called the Red Set.

In the present work a set $\{\mathscr R^{\psi(\mathcal C)}\}$ of \ ``Reverse Tableaux'' is constructed through an Enabling Proposition.  They are shown to define the same components as the component tableaux and furthermore
give a factorisation of each new invariant in a chosen sequence, via successive linearisation of preceding invariants.  Via Krull's theorem this factorisation is to give the required surjectivity.

There can be several reverse tableaux for a given Red Set, yet each determine the same variety as the component tableaux with the given Red Set and define the same components. The flexibility of having several ``equivalent'' reverse tableaux, absent from the more ``rigid'' component tableaux, is quite essential for factorisation.

This procedure is quite different from the classical approach to surjectivity requiring a \textit{geometric} description of the nilcone, which in the present case seems unattainable.

Nothing of this complexity has ever been tackled before and indeed the methods used here are entirely new.

\

     \textbf{Acknowledgements}.

     \

     The work of the first author was supported by ISF grant 1957/21 jointly held by M. Gorelik, Weizmann Institute and S. Reif, Bar-Ilan University, by ISF grant No. 1347/23 held D. Novikov at the Weizmann Institute and by a grant from the Mathematics Faculty, Weizmann Institute.

     \

     The main results described here were the subject of an invited talk by Y. Fittouhi during December 2024 and invited talks by A. Joseph  in the Weizmann Institute during December 2024, April 2025, and by Y. Fittouhi in Bar-Ilan University during November 2025.

\section {Introduction} \label {1}

\subsection {Avant-propos}\label {1.1}

To paraphrase Leonardo da Ponte, a component is like the Phoenix, everyone swears it exists but nobody has ever seen one.

Indeed components of a given algebraic variety, exist via the Lasker-Noether theorem \cite [Thm. 4.5]{AM} but in reality are notoriously difficult to find. Their search is a far-reaching generalisation of finding the prime factors of a positive integer, an age-old problem now viewed of added importance for its role in breaking public key encryption.

A further important related development occurs in the context of certain non-commutative algebras, notably enveloping algebras and their quantum analogues.  In this vast amounts of work has been published on prime (and primitive) ideals and the decompositions to which they lead.  This is partly because of their importance in Representation Theory and partly because they are fewer so more amenable to classification. In this P. Gabriel et al \cite {BGR} and J. Dixmier \cite {D}, made notable early  contributions and expositions.

Here we consider the action of a parabolic subgroup $P$ on the Lie algebra $\mathfrak m$ of its nilradical. The relevant invariant ring is $\mathcal I:=\mathbb C[\mathfrak m]^{P'}$, with $P'$ being the derived group of $P$.  Via Richardson's theorem it is polynomial \cite [2.2.3]{FJ1}.

We extend our study \cite {FJ5} in type $A$ of the components of the ``nilfibre'' $\mathscr N$, defined to be the zero locus in $\mathfrak m$ of the augmentation $\mathcal I_+$.

In type $A$, Benlolo and Sanderson \cite {BS} suggested a choice of polynomial generators for $\mathcal I$.  This was verified in \cite {JM}.  We call them the Benlolo-Sanderson (BS) invariants.

In \cite {FJ5} we constructed a map, called the component map, from a set of explicit combinatorial data to the set of components of $\mathscr N$.  It was shown to be injective \cite [Sect. 7]{FJ5}.

We propose to eventually show that the component map is also surjective.  In the resulting endeavour, we simplified the said combinatorial data as a multi-set, called simply the ``Red Set'', which assuming surjectivity gives a more transparent classification of components of $\mathscr N$ in type $A$.

To have any chance of success of finding \textit{all} components of a variety, conventional wisdom would say that  one needs a \textit{geometric} understanding of that variety and not just its description as the zero locus of an algebraic set. Ideally the variety should admit an algebraic group action so that the components are the closures of its orbits.  This is the case of the Springer fibre for which such a description results from the Steinberg triple variety.  Thereby its components were  determined \cite {S},\cite {St} with relative ease.

For $\mathscr N$ we have so far no convenient geometric description as in the case of the Springer fibre, and little hope of one.   Even in type $A$, the nilfibre can be far from irreducible, unlike the case, studied in the classical work of Kostant \cite {K}, when $P$ is the entire simple algebraic group acting on its Lie algebra.  Again for the Springer fibre, the components were, under a simple transformation, Lagrangian subvarieties of associated $G$ orbits, whilst in the present case, one component of the same nilfibre might be Lagrangian whilst another need not be, yet both have the same dimension  \cite [8.2]{FJ5}.

Thus it is hardly surprising that no-one has even dared to determine the components of $\mathscr N$ before.  Such an endeavour could have been initiated even one hundred years ago.

Yet as we saw the problem simplifies through Richardson's theorem \cite [Sect.2]{FJ1}, going back ``only'' fifty years.


Our present work might well be read in another fifty years time for those wishing to extend the present results to general type, which we anticipate will be a huge challenge and a fascinating endeavour.

Already the determination of the invariant generators in classical type constituted a remarkable PhD thesis of Elena Perelman \cite {Pe}, which has not been reached let alone bettered over the last twenty-three years.  Besides generators are barely understood particularly for general type \cite [3.3]{FJ1}.

This work is a continuation of \cite {FJ5}, having itself origins in \cite {FJ1}, \cite {FJ2}, \cite {FJ3}. These earlier projects were to describe Weierstrass sections for the action of a parabolic subgroup of $SL(n)$ on the Lie algebra of its nilradical and in which a ``canonical component'' tableau was introduced \cite {FJ4}.  Their main interest to the present reader lies in the fact that they motivated the present ``totally out of blue construction'' described in \cite {FJ5} and significantly deepened here, through the reverse tableaux, which were an inspired invention by the first named author.

\subsection {Notation and Conventions}\label {1.2}

These will follow closely those of \cite {FJ5} but will be redefined here as necessary.  For the moment we note that:

The base field is assumed to be algebraically closed and of characteristic zero.  It will simply be denoted as $\mathbb C$.   Given $m\leq n$ is  positive integers, set $[m,n]:=\{m, m+1,\ldots,n\}$.  If $m$, $n$, or both are omitted we note this respectively as $]m,n],[m,n[,]m,n[$.

 \section {The Main Players}\label {2}

 \subsection {Matrices}\label {2.1}

 Recall the notation of \cite [Sections $2,3$]{FJ5}. Identify $\mathfrak {gl}(n)$ with the space $\textbf{M}_n$, or simply \textbf{M}, of $n \times n $ matrices under commutation, and $\mathfrak {sl}(n)$ its Lie subalgebra of trace zero matrices.

Let $\mathfrak b$ be the Borel subalgebra of $\mathfrak {gl}(n)$ of upper triangular matrices and let $\mathfrak p \supset \mathfrak b$ be a (standard) parabolic subalgebra of $\mathfrak {gl}(n)$.  The latter is completely determined by a composition  $(c_1,c_2,\ldots,c_k)$ of $n$ defining the Levi factor $\mathfrak r$ of $\mathfrak p$ as the vector space spanned by the co-ordinate vectors in the $\{c_i\times c_i\}_{i=1}^k$ blocks $\textbf{B}_i$ down the diagonal of $\textbf{M}_n$.  The nilradical $\mathfrak m$ of $\mathfrak p$ is the space spanned by the upper triangular part of $\textbf{M}_n$ strictly above these blocks.

%
%

\subsection {Diagrams, Columns and Rows}\label {2.2}

Define a diagram $\mathscr D$ with columns \newline $C_1,C_2,\ldots,C_k$ satisfying $\height C_i=c_i$.

 For each column $C_i$, let $\textbf{C}_i$ denote the rectangular block in $\mathfrak m$ lying above $\textbf{B}_i$.  We call it the $i^{th}$ column block. Its width is $c_i$ and its height is $\sum_{j < i}c_j$.  In particular $\textbf{C}_1=0$.

    Thus  $\mathfrak m:= \oplus_{i=2}^k \textbf{C}_i$.

    Let $x_{i,j}$ denote the $ij^{th}$ entry of \textbf{M}.  They are called co-ordinate vectors and form its standard basis.  Let $x^*_{i,j}$ denote the function taking value $1$ on $x_{i,j}$ and $0$ on the remaining co-ordinate vectors.

Let $R_i:i=1,2,\ldots,$ denote the rows $\mathscr D$.  We view them as \textit{descending down} $\mathscr D$ as $i$ increases. Set $R^u=\cup_{i=1}^u R_i$ and $R^{>u}$ the set of remaining rows, that is not in $R^u$.  Set $b_{i,j}:=C_i \cap R_j$, called the $ij^{th}$ box of $\mathscr D$.

\subsection {Adjacent and Neighbouring Columns} \label {2.3}

  Successive columns of $\mathscr D$ are said to be adjacent. Two boxes in the same row in not necessarily adjacent columns of $\mathscr D$, with no boxes between them, are called adjacent.

  If $C,C'$ are columns of $\mathscr D$ we shall always mean $C$ to lie to the left of $C'$.  Denote by $[C,C']$ the set of columns of $\mathscr D$ between $C,C'$. If $C$, $C'$, or both are omitted, then we denote the resulting set as $]C,C'],[C,C'[,]C,C'[$.  According to context we may also mean this to be the set of boxes in these columns or their entries (see \ref {2.5}).

   Two columns $C,C'$ of $\mathscr D$ having height $s$  are said to be neighbouring \cite [4.1.2]{FJ1}, if they are no columns of height $s$ strictly between them, that is in $]C,C'[$.  Let $\textbf{g}$ be the number of such neighbouring pairs.

   The lower part $C^{>u}$ of a column $C$ relative to $u \in \mathscr N$ is defined to be $R^{>u}\cap C$.



\subsection {Rectangles} \label {2.4}


 If $C,C'$ are neighbouring columns of height $s$, we set $R^s_{C,C'}=R^s\cap [C,C'] $, sometimes referred to as the rectangle defined by these columns.  Define the left rectangle to be $^lR^s_{C,C'}:=R^s_{C,C'}\setminus C=R^s\cap]C,C']$.

 \subsection {Tableaux}\label {2.5}

Define a tableau $\mathscr T$ by filling $\mathscr D$ with distinct entries from $[1,n]$ going down columns and then from left to right (cf \cite [2.2.2]{FJ5}).  Observe \cite [2.2.2]{FJ5} that in this, right going lines $\ell_{i,j}$ with entries $i,j$ in distinct columns of $\mathscr T$ define a basis  of $\mathfrak m$ formed from the $x_{i,j}$.  If there is no such restriction on $i,j$ these elements form a basis of \textbf{M}.

Observe that $x_{i,j}$ is the root vector of the root $\alpha_{i,j}$ both being given by the line $\ell_{i,j}:i\neq j$. To simply our terminology we shall conflate a root vector with a root.

 \subsection {Surrounding columns}\label {2.6}

 Consider $C_{r+1}:r \in [1,k-1]$. A pair of neighbouring columns $C,C'$ is said to surround $C_{r+1}$ if $C_{r+1} \in ]C,C']$.  This implies that $C_r \in [C,C'[$ and so $C,C'$ surround the pair $C_r,C_{r+1}$ in the sense of \cite [3.1.3]{FJ5}.

 \begin {lemma} There exists at most one pair of neighbouring columns of a given height, say $s$, which surround $C_{r+1}$.
 \end {lemma}

 \begin {proof} This is an obvious consequence of the fact that a pair of neighbouring columns of a given height $s$ cannot overlap except by sharing a single common column.
 \end {proof}

 \subsection {The Benlolo-Sanderson Invariants}{\label {2.7}

 A Benlolo-Sanderson (BS) invariant $I^s_{C,C'}$ is given by a pair $C,C')$ of neighbouring columns of height $s$.  They are essentially determinants and in particular multi-linear.  They were first defined in \cite {BS} and further studied in \cite {JM}.

  The (BS) invariants are (a choice of) generators of $\mathcal I:=\mathbb C[\mathfrak m]^{P'}$, which is polynomial on $\mathbf {g}$ generators.

  Their common zero locus is (defined to be) $\mathscr N$.

\begin {lemma}  The number $d(R^s\cap ]C,C'])$ of entries in $^lR^s\cap ]C,C']$ is just the degree $d$ of $I^s_{C,C'}$, viewed as a polynomial in the $x_{i,j}$.
\end {lemma}

\begin {proof}

   As noted in \cite [Lemma 4.2.5]{FJ1} (further detailed \cite [Prop. 5.1.1(i)]{FJ5} and \cite [Lemma 4.3.2]{FJ5}) $I^s_{C,C'} $ is homogeneous and a sum of monomials one of which one is obtained by taking all disjoint composite lines in $R^s_{C,C'}$ joining the boxes of $C$ to those of $C'$, going from left to right and passing through all the entries of $R^s\cap [C,C']$, once and only once.

 This is just a simple consequence of the rules to compute a determinant as noted in the footnote to \cite [1.3]{FJ4} and since only one monomial is selected even signs do not matter


  The disjoint composite lines in the paragraph preceding the lemma, are formed from $|[C,C']\setminus C|$ strictly right going lines giving a monomial of degree $d(R^s\cap ]C,C'])$.  Hence the assertion.
\end {proof}

  \subsection {Nachbarschaftsanschauung}{\label {2.8}

  With hindsight and tongue-in-cheek we can formulate our ``neighbourhood philosophy''.

  It is that the components of $\mathscr N$ are obtained from $\mathscr T$ to give ``component tableaux'' $\mathscr T^\mathcal C$ by repeating numerical entries which first appear in $\mathscr T$ on the right and descending by rows specified by surrounding pairs of neighbouring columns.

  Logical - as the generators of $\mathcal I$ are given through neighbouring columns.  Yet it is  not so obviously discoverable even by ``artificial intelligence''.  Our own success, documented in \cite {FJ4}, resulted from ``serendipity''. It makes a fundamental contribution to the present endeavour.

  Again reverse tableaux are also constructed through neighbouring columns but by a different procedure.

 \section {Component Tableaux}\label {3}

 \subsection {Line Decoration}\label {3.1}

 A component tableau $\mathscr T^\mathcal C$ is defined by combinatorial data $\mathcal C$ given in \cite [3.2]{FJ5}.  It is viewed as the tableau $\mathscr T$ decorated by lines $\ell_{i,j}$, joining entries $i,j$ in boxes of distinct columns and going  from left to right.

 The lines $\ell_{i,j}$ may be labelled with a $1$, a $\ast$ or not at all.  The subset with label $1$ (resp. $\ast$) is denoted by $S(\mathscr T^\mathcal C)$ (resp.  $Y(\mathscr T^\mathcal C)$).

 Then evaluation of \textit {any}  Benlolo-Sanderson invariant through this labelling replaces it by one of the co-ordinates with label $\ast$, in other words it gives  a Weierstrass section \cite [Thm. 5.2.6]{FJ5} for the component defined by the component tableau $\mathscr T^\mathcal C$.

 Such a labelling was found in \cite {FJ1} and was a remarkable achievement of the first named author.  Notably the time-honoured Kostant construction, developed for the ``Toda lattice'', fails miserably in the present context \cite [1.3]{FJ5}.  The work in \cite {FJ1} used all sorts of strange gizmos, yet was finally understood in \cite {FJ4} by what can be now be incorporated into our neighbourhood philosophy \ref {2.8}.

 Then a programme was lanced to generalize this construction. It culminated in the rules in \cite [3.2.2, 3.2.6, 3.2.7]{FJ5}.  These are explained below.

As a by-product of the construction of Reverse Tableaux the (mysterious) combinatorial data encoded in $\mathcal C$ was reformulated in terms of a much more transparent multi-set called simply the ``Red Set'', \ref {3.2.4}.


\subsection {The Component Map}\label {3.2}

Associated to $\mathscr T^\mathcal C$ there is a set $X(\mathscr T^\mathcal C)$ of excluded root vectors in $\mathfrak m$, described explicitly in \cite [4.3.9, 4.3.10]{FJ5} and extended in \ref {4.4}. Let $\mathfrak u^\mathcal C$ be the subspace of $\mathfrak m$ spanned by vectors which are \textit{not} excluded and set $\mathscr C:=\overline {B.\mathfrak u^\mathcal C}$.  By \cite [Sect.6]{FJ5}, $\mathscr C)$ is a component of the nilfibre $\mathscr N$.   The map $\mathscr T^\mathcal C\mapsto \mathscr C$, is called the component map.  It is injective \cite [Thm. 7.5]{FJ5}.  Proving it to be surjective will complete our description of the components of $\mathscr N$, making the Red Set the classifier of components.

We believe our construction to be straightforward if not natural, but not having been exposed to it during several years some readers may flounder, become submerged and drown in the details.  Therefore we shall attempt a ``Royal Route'' to the construction replete with illustrated step by step instructions diagrammically presented in Section \ref {10}.

\subsubsection {The Construction of the Component Tableaux \cite [3.2.1]{FJ5}}\label {3.2.1}

The construction of $\mathscr T^\mathcal C$ involves an intermediate step through a semi-standard component tableau $\mathscr T^\mathcal C(\infty)$ produced stepwise \cite [3.1.4]{FJ5}. In this, $\mathscr T$ and $\mathscr T^\mathcal C$ are viewed as sub-tableaux of $\mathscr T^\mathcal C(\infty)$ with the first step in the construction of $\mathscr T^\mathcal C(\infty)$  being $\mathscr T$ itself.

We regard $\mathscr T^\mathcal C(\infty)$ as evolving in time (the latter having only a positive direction).  Inspired \cite[3.1.2]{FJ5} by the \textit{vav conversive} of ancient texts, one may cut down notation by not slavishly specifying the exact point in time reached.  The reader may appreciate this convention, some may find it insufficient.  To avoid any ambiguity we did not use this convention in \cite [3.2.1-3.2.5]{FJ5}.

\subsubsection {Batches \cite [3.1.4]{FJ5}}\label {3.2.2}

As $\mathscr T^\mathcal C(\infty)$ evolves in time, its empty boxes fill up in each row $R_s$, with then $s$ increasing, following the procedure described below.  A particularly simple case is described in \cite [4.7]{FJ4} which the reader might find helpful to study.  It is referred to as the canonical tableau and defines the ``canonical component''.

Let $C,C'$ be neighbouring columns of height $s$.  The batch $\mathscr B^s_{C,C'}$ is defined to consist of the rightmost entries having a given value $j$ in $R_s\cap [C,C'[$, after the row $R_s$ is filled up.


\subsubsection {Strings}\label {3.2.3}

Each $i \in [1,n]$ moves rightwards through successive columns of $\mathscr T^\mathcal C(\infty)$ from its starting point in $\mathscr T$, into empty boxes.  They are deemed to form an ``$i$-string'', denoted $S(i)$ (\cite [Definition 4.2.1]{FJ5}). Conversely very entry $j$ in $\mathscr T^\mathcal C(\infty)$ belongs to a \textit{single} $j$-string. A singleton $i$ is deemed an $i$-string, but in general there can appear several entries of $i$ in $\mathscr T^\mathcal C(\infty)$.  As $i$ strings move rightwards from their origin in $\mathscr T$, the unique leftmost element of an $i$-string lies in $\mathscr T$.

A line joining two boxes in adjacent columns having the same numerical entry is called a neutral line \cite [5.2.1]{FJ5}.  These lines may be viewed as ``tracing'' the $i$ strings.

Given $i_1,i_2 \in [1,n]$ distinct, then $S(i_1),S(i_2)$ cannot cross \cite [Lemma 4.2.1]{FJ5}.  If they both pass through a common column with $S(i_2)$ below $S(i_1)$, necessarily strictly, then $S(i_2)$ starts in a column to the left of $S(i_1)$ \cite [4.2.2]{FJ5}.

An $i$-string may move horizontally, or down by $m\geq 1$ rows on passing from a column $C_r$ to the next $C_{r+1}$ given by specific rules \cite [3.2]{FJ5}, using pairs of neighbouring columns.


In this a pair of neighbouring columns cannot be used twice and a pair not yet used in this step-wise procedure is called \textbf{free} \cite [Definition 3.1.5]{FJ5}. Eventually all such pairs must be used \textit{exactly once}. In particular one and only one element is chosen from each batch $\mathscr B^s_{C,C'}$.  Within a given row (that is for fixed $s$) the order of choice is immaterial; but the order  of going down rows must be respected \cite [3.1.5]{FJ5}. For one thing $\mathscr B^s_{C,C'}$ depends on choices of elements from batches made in previous rows.


\

Let us expand thie above description by explaining the consequences of \cite [3.2.2, Rule (1)(a)]{FJ5}.

Suppose $i\in R_s\cap C_r$. At this point in time all columns to the right of $C_r$ have height $\geq s-1$. Let $R_{s+m}\cap C_{r+1}$ be the first empty box of $C_{r+1}$, going \newline downwards.

Suppose $m=0$. Then $i$ is moved into the empty box $R_{s}\cap C_{r+1}$.

Suppose $m > 0$.  Take $m' \in \mathbb N$.

Suppose there are $m+m'$ pairs of free neighbouring pairs $C,C'$ of heights exactly  $s+v:v \in [0,m+m'-1]$ which surround $C_r,C_{r+1}$ in the sense of \cite [3.1.3]{FJ5}, that is to say both $C_r,C_{r+1}$ lie in $[C,C']$.

Then (\cite [3.2.5]{FJ5}), there can be at first a downward step by $m\geq 1$ rows from a column $C_r$ to its right adjacent column $C_{r+1}$ into which $i$ enters a first empty box.  This leads to $m$ vertical lines $\ell_{i,j}$, labelled by a $\ast$  going from $i$ in this newly filled box going upwards to the set of necessarily distinct entries in the subset $J^\mathcal C_{r+1}\subset C_{r+1}$ of the lowest $m:=|J^\mathcal C_{r+1}|$ boxes of $C_{r+1}$.  (Hence in the obvious sense $J^\mathcal C_{r+1}$ is a connected subset of a lower part of $C_{r+1}$.)

After that (\cite [3.2.5]{FJ5}), there can be a further $m'-1$ downward steps of \textit{just one row} on going from $C_r$ to the first empty box in $C_{r+1}$ and \textit{from then on} no steps at all.  Then there are $m'$ vertical lines $\ell_{i,j}$  labelled by a $\ast$ each from one of these new boxes to the lowest box in $C_{r+1}$ (giving it multiplicity $m'$).

Corresponding to each line $\ell_{i,j}$, there is a co-ordinate vector $x_{i,j}\in \textbf {M}$ labelled by a $\ast$.

If is no collection  of $m$ pairs of free neighbouring  pairs of heights $s+v:v \in [0,m-1]$, then the $i$-string is ``stopped'' at $C_r$ - \cite [3.2.2,Rule(2)]{FJ5}. In this case, $i$ is joined by a line $\ell_{i,j}$ with label $1$ to the unique element lowest $j \in C_{r+1}$, not already joined by a right and necessarily up-going line with label $1$, \cite [3.2.1, 3.2.6]{FJ5}, defined at a previous step.

Correspondingly the co-ordinate vector $x_{i,j}\in \textbf {M}$ is labelled by a $1$.

These rules derive automatically from the condition that a given pair of neighbouring columns cannot be used twice.  At first the reason for this might not seem obvious but it does indeed become increasingly apparent - stay with us!

The resulting configuration of lines is illustrated in \cite [Sect. 9, Fig. 1]{FJ5}.

In the above there can be several possibilities defined by choices in the batches,  Each lead to a component of $\mathscr N$ defined by the excluded roots (\cite [4.2.3, 4.2.4]{FJ5}, \ref {4.3.6}) and the above labelling by $1,\ast$ gives rise to a Weierstrass section \cite [Thm. 5.2.6]{FJ5}, associated to that component.

In this, the excluded roots are encircled in \textbf{M}.  All the root vectors labelled by $\ast$ are encircled and none labelled by a $1$ \cite [Lemma 4.1.5. Prop. 4.2.5]{FJ5}.  This is recalled in \ref {8.1}.

An example comes from \cite [Sect. 9, Example 5]{FJ5}, which comes from the composition $(3,2,1,2,2,1,3)$.  In this between $(C_4,C_5)$ the entry $7$ goes down by two rows using the pairs $(C_3,C_6),(C_3,C_4)$ of heights $1,2$ and then $5$ goes down by one row using $(C_1,C_7)$, of height $3$.  However instead of this last step, $8$ could not then have gone down by two rows as this would require using $(C_3,C_4)$ again (as well as using $(C_1,C_7)$).

The examples in \cite [Sect. 9]{FJ5} illustrate the decoration in $\mathscr T$ given by $\mathcal C$ and the resulting labelling in \textbf{M}.

Note in particular that an entry of \textbf{M} cannot both be labelled a $\ast$ and a $1$.   This rule is \cite [Cor. 3.2.6(i)]{FJ5}.  It is illustrated in \cite [Sect. 9, Fig. 1]{FJ5}, in which one may observe that a line with label $1$ \textit{cannot} have the same end-points as the composition of neutral line and a red line - taking care in this to note the fourth sentence of the caption to  \cite [9, Fig.1]{FJ5}. (In this we have perhaps been guilty of trying to cram too much information into this one figure!)

\subsubsection {The Red Set}\label {3.2.4}

The elements of $J_{r+1}^\mathcal C:r \in [1,k-1]$ are coloured in red and the remaining elements of $\mathscr T^\mathcal C(\infty)$ in black (simply because these were the colours most visible and the pens readily available.


  $(*)$. Consequently for a given numerical label $i$, there is at most one entry in red and at least one entry in black, possibly several.   Notice in this, the black entries with a given numerical value (say $i$) are joined by neutral lines to form an $i$-string.

\

\textbf{Example $1$}.  Consider the composition $(1,2,3,1,1,3,2)$.  We compute the six Red Sets for the component tableaux through a flow chart in Sect. $10$, Figure $1$.
\


  The reader is advised to pencil in the neutral lines and those with label $1,\ast$ respectively in black and red ink in these examples.  This is not just fun but helps to follow the arguments of this paper.  It is the first step on the Royal Route to tableaux.

  \

  The cardinality of $J^\mathcal C_{r+1}$ as a \textit{set} is called its multiplicity $m^\mathcal C_{r+1}$.  The superscript is often omitted.  The subscript is sometimes omitted.  We need only be concerned with the case when $J^\mathcal C_{r+1}$ is not empty, so then $m\geq 1$.

An $i$-string is completely determined by the entry of $i \in \mathscr T$, by these downwards steps and the place where it is stopped (if at all).

Given that $J^\mathcal C_{r+1}$ is non-empty, then by \cite [3.2.5]{FJ5} at most its lowest entry $j^\mathcal C_{r+1}$
can lie at the top of strictly more than one $i$-strings meeting $C_{r+1}$ - see \cite [Sect. 8, Fig. 1]{FJ5}.

Again this property can be traced to the condition that a pair of neighbouring columns can only be used once.

We call the multiplicity $m'_{r+1}$ of $j_{r+1}$. It can be simply equal to  $1$. The subscript is sometimes omitted.  Note that $m_{r+1}'\geq 1$, if $J_{r+1}\neq \phi$.

Set $h_{r+1}=\height C_{r+1}$.  Again the subscript is sometimes omitted.

Let $R(\mathscr T^\mathcal C)$ denote the multi-set which is the union of the multi-sets  $J_{r+1}^{\mathcal C}$, for $r \in [1,k-1]$ of $\mathscr T$.  These entries will be coloured in red and called the Red Set of $\mathscr T^\mathcal C(\infty)$.  Here capital letters are used to emphasize that this is a definition more than a description.  Indeed the Red Set is a multi-set but we believe that the more precise term ``red multi-set'' just too clumsy. Yet multiplicities \textit{must} be taken into account for Lemma \ref {3.3} to be valid.

Finally $\mathscr T^\mathcal C$ can be recovered from $\mathscr T^\mathcal C(\infty)$ by collapsing $i$-strings \cite [5.2.1]{FJ5} and retaining the numerical labels.

\

In the right hand column of Figure $10$ of Sect. $10$, the reverse tableaux are listed for the component tableaux of Figure $1$ of Sect. $10$.  Those with the same Red Set can differ; but by Theorem \ref {9.3} all define the same component, for the corresponding component tableau.

\subsubsection {The Red Lines}\label {3.2.5}

\

Consider a vertical line from $C_{r+1}(\infty):r \in [1,k-1]$ to the red subset $J^\mathcal C_{r+1}$.  It takes a black entry of $i\in C_{r+1}(\infty)\setminus C_{r+1}$, to a red entry of $j\in J^\mathcal C_{r+1} \subset C_{r+1}$ and is labelled by a $\ast$.  Both end-points are determined by their numerical values.

The red line $\ell^\mathcal C_{i,j}$ with label $\ast$ joining the unique entries in $\mathscr T^\mathcal C$ is obtained by collapsing  the  neutral lines between this $i\in C_{r+1}(\infty)\setminus C_{r+1}$ and the unique leftmost copy of $i$, which lies in $\mathscr T^\mathcal C$.  Since in this $i$ lies in a column strictly to the left of $C_{r+1}$ which contains $j$, we obtain $i<j$.

This leads to the following

\

\textbf{Definition}.  A red line $\ell^\mathcal C_{i,j}$ in $\mathscr T^\mathcal C(\infty)$ joins the unique red entry of $\mathscr T^\mathcal C(\infty)$ with numerical value $j$ to the unique leftmost black entry with numerical value $i$.

\

It determines a unique line $\ell^\mathcal C_{i,j}$ joining entries $i<j$ of $\mathscr T^\mathcal C$ with label $\ast$.

\

By \cite [3.2.5]{FJ5} the number of red lines in $\mathscr T^\mathcal C$ is the number of pairs of neighbouring columns and so the number $\textbf{g}$ of generating invariants.



\subsection {The Red Map}\label {3.3}

The map $ \mathscr T^\mathcal C \mapsto R(\mathscr T^\mathcal C)$ is called the red map.

\begin {lemma} The red map is injective.  In other words $ R(\mathscr T^\mathcal C)$ uniquely determines $\mathcal C$.
\end {lemma}

\begin {proof} Each $i$-string encodes the batches to which $i$ belongs.  Indeed $i\in R_s$ descends to $R_{s+l}$ as it passes across adjacent columns $C_r,C_{r+1}$ of $\mathscr T^\mathcal C(\infty)$ exactly when $i$ is chosen from the batches $\mathscr B^{s+t}_{C_{k_t},C'_{k_t}}$ given by $m$ free pairs of neighbouring columns $C_{k_t},C'_{k_t}$  of height $s+t:t \in [0,m-1]$, surrounding $C_{r},C_{r+1}$ \cite [3.1.3,3.2.3]{FJ5}.  Thus the set of descents of the $i$-strings determines $\mathcal C$.

On the other hand, viewed as a set $J^\mathcal C_{r+1}$, determines which $i$-strings descend by $m$ rows as they pass between $C_{r},C_{r+1}$ entering the first empty box below $C_{r+1}$. More precisely $i\in R_{h-m+1}\cap C_r(\infty)$ descends by $m$ rows into $R_{h+1}\cap C_{r+1}(\infty)$.

Secondly from free neighbouring columns of heights $h+t-1:t \in [1,m']$ surrounding $C_r,C_{r+1}$, there are further $m'$ generally different $i$-strings for which $i$ is the entry of $R_{h+t-1}\cap C_r(\infty)$, descending each by one row into $R_{h+t}\cap C_{r+1}(\infty)$.
It is exactly the previously defined $i$-string when $t=1$.  (Thus if $m\geq 1$, then $m'\geq 1$.)

In brief $ R(\mathscr T^\mathcal C)$, with \textit{multiplicities counted!} determines the descents of each $i$-string which in turn determines the batches to which $i$ belongs and this defines $\mathcal C$.

Hence the assertion of the lemma.

\end {proof}

\textbf{Remark 1.} Altogether this description requires neighbouring columns of heights $i:i \in [h-m+1,h+m'-1]; h =\height C_{r+1}$, surrounding $C_r,C_{r+1}$, which are ``free'', that is to say not used at a previous step \cite [3.1.5]{FJ5}. Recall that we may assume $m,m'\geq 1$ and so there are no gaps in this set of heights!

\

This description emphasizes the proof in \cite [3.2.5]{FJ5} that only $j^\mathcal C_{r+1}$ amongst the boxes of $J^\mathcal C_{r+1}$ can have multiplicity $>1$. Otherwise we would be using twice some pair of neighbouring columns, as the reader may confirm on gentle reflection. As we shall see this phenomenon also arises for the reverse tableaux. Again (\ref {3.2.3}), this is the price of freedom!

\

\textbf{Remark 2.} The descent of the $i$-strings referred to in the lemma is illustrated graphically by \cite [Sect. 9, Fig. 1]{FJ5}.  The $i$-string descending by $m=m_2$ rows is shown by the dashed green line starting from the highest red box. (Again ``attend and mark'' the fourth sentence of the caption.) The remaining dashed green lines (except the very lowest) descend by one row.

\

\textbf{Remark 3}.  One might like to determine all possible choices of Red Sets with respect to a given composition.  This is a matter of exactly distributing all pairs of neighbouring columns so that each pair contributes exactly once.  In the present work we do not need to address this knotty problem, which is nevertheless solved in \cite {FJ7}.

\

\textbf{Remark 4}. We recommend that the reader become familiar with the algorithm determining the Red Set.  We ourselves did altogether several hundred cases.

\

Some Red Sets for component tableaux are given in Section \ref {10}, Example $1$.


%

\

%

 \textbf{Remark 5}.  One may say that the Red Set gives an alternative way to incorporate the combinatorial data determining the component tableaux.  Yet  the reader may need to return to the batches to determine the Red Set and slavishly follow the analysis in \cite [3.2]{FJ5}, otherwise it is not obvious how to determine $S(\mathscr T^\mathcal C), Y(\mathscr T^\mathcal C)$, even though they are implicit.
  The examples in \cite [Sect. 9]{FJ5} might help, but here we were rather brief. Like piano you have to practice\footnote {As Miss Elizabeth Bennett pointed out to Fitzwilliam Darcy.}.

 \subsection {A Total Order Relation on the Red Set}\label {3.4}

 We define an order relation on a Red Set by moving in $\mathscr T^\mathcal C$ from left to right and then \textit{upwards} within a given column.  This is nearly but not quite the numerographical order with respect to the numbering of boxes in $\mathscr T$.  Moving upwards is natural when constructing a reverse tableau since by \ref {4.3}$(*)$ we must first apply its construction to the lowest element black entry in a given column, which is then recoloured in red. In addition $j^\mathcal C_{r+1}$ is linearly ordered as described below.

 Recalling that $i$-strings do not cross (\ref {3.2.3}), those passing through $C_{r+1}$ are linearly order by the  multi-set $j^\mathcal C_{r+1}$.  Those lower down start further to the left and if strictly their entries decrease, as in Example $2_1$. Otherwise they increase. Thus their numerical entries may oscillate, as in the following example illustrates.  Yet we may reassure the reader that this funny behaviour has little influence!

 \

 \textbf{Examples $2_1$}.  Consider the composition $(1,2,1,2)$.   Using the neighbouring pair $(C_1,C_3)$, one may lower $2$ into $R_2\cap C_3$ and then through the neighbouring pair $(C_3,C_4)$, further lower $2$ into $R_3\cap C_4$. This gives the Red Set $(4,6)$.

Again there is a vertical line from the second $2$ to $4$ and from the third $2$ to $6$.  Then $3$ is stopped by $4$, whilst $1$ is stopped by $2$ and $4$ by $5$ which after collapsing neutral lines giving lines $\ell_{1,2},\ell_{3,4}, \ell_{4,5}$ with label $1$ and lines $\ell_{2,4}, \ell_{2,6}$ with label $\ast$.

Anticipating covering (\ref {9.2.1}), we note that in this case $x_{4,6}$ is an excluded root covered by $x_{4,5}$.

In the second case the Red Set is $(4,4)$.  The excluded root $x_{1,4}$ is covered by the root $x_{1,2}$ with label $1$.

 Consider the composition $(1,2,1,2)$.  There is a component tableau $\mathscr T^\mathcal C$ in which $2$ goes below $4$ and $3$ below this $2$ in $C_3(\infty)$.  In this $j^\mathcal C_3$ is the (red) multi-set $(4,4)$, whilst in the reverse tableau, a red $4$ appears twice.

On the other hand, in the composition $(2,1,1,2)$  with Red Set $(4,4)$, the entries $2,3$ in $C_3(\infty)$ appear in reverse order.

Consider the composition $(2,3,1,1,3,2)$ with Red Set $(7,7,7)$,  Going downwards $C_4(\infty)$ has entries $6,4,5$ in that order.


 \subsection {Complete Sequences of Pairs of Neighbouring Columns.}\label {3.5}
Fix a diagram $\mathscr D$ and let $\mathscr P$ (resp. $\mathscr P_t$) denote the set of pairs of neighbouring columns (resp. of height $t$).  Set $ \textit{p}=|\mathscr P|$.

Label the columns of $\mathscr D$ as in \cite [2.1.1]{FJ5}.  Then a particular element of $\mathscr P_{s_i}:i\in \textit{p}$, may be denoted by $(C_{r_i}(s_i),C_{r'_i}(s_i))$ with $r_i, r'_i:r_i<r'_i \in [1,k]$.  Here the common heights $s_i$ of a given pair, will sometimes be omitted from the notation.

\

\textbf{Definition.} A complete sequence of pairs of neighbouring columns is a total order on $\mathscr P$ viewed as a sequence $\{C_{r_i},C_{r'_i}\}_{i=1,2,\ldots, \textit{p}}$ of various heights.

\

The construction of a reverse tableau will follow such a complete sequence.  These sequences can be different for different Red Sets of the component tableaux and to ensure that all reverse tableaux corresponding to the various component tableaux, all may be called upon.

\

The first example in $2_1$ gives the Red Set $(4,4)$ coming from the complete sequence $(C_1,C_3);(C_2,C_4)$, whilst it also gives the Red Set $(4,6)$, so only one complete sequence needed.

\

We had optimistically thought that only one sequence was needed in all cases following some natural total order on $\mathscr P$. However this is insufficient to obtain the set of all reverse tableaux.  The simplest example we found was for the composition $(3,2,1,1,3,2)$.  This has four components. A different one is missed using either the sequence $(C_3,C_4),(C_1,C_5),(C_2,C_6)$ or the sequence $(C_3,C_4),(C_2,C_6),(C_1,C_5)$.  \textit {In general a suitable union of sequences derives from the total order on the Red Set}, we show this for the reverse tableaux constructed here.  The indefatigable reader may wish to check this in this last example.   A slightly more complicated example is given in Example $10$ of Section \ref {10}.

Recall that $\mathcal I=\mathbb C[\mathfrak m]^{P'}$ is a polynomial algebra and its generators can be assumed of the form $I^s_{C,C'}$ where $C,C'$ are neighbouring columns of height $s$.

 In particular a total order on the Red Set induces a total order on $\mathscr P$ which in turn induces a total order on the (BS) generators $\{I_1,I_2,\ldots, I_\textbf{g}\}$.

 \

The surjectivity of the component map (\ref {3.2}) would mean that the Red Set is the classifying set for components of the nilfibre $\mathscr N$.

 \subsection {The Elephant in the Room.}\label {3.6}

  Given a total order on generators assume that we have identified a component $\mathcal C$ of the zero locus of $\mathscr Z_{i-1}$ of the first $i-1$ generators from this list and consider the Weierstrass section it defines.  Here we may replace each $I_j:j\in [1,i-1]$ by a linear function  $\ell(I_j)$ defined by this Weierstrass section. Then by Krull's theorem the factors of the $i^{th}$ invariant gives the components of the variety of zeros of $I_i$ restricted to $\ell(I_j)=0:j=1,2,\ldots, i-1$.

  Now from each such factor (and the components of $\mathscr Z_{i-1}$) assume that we have identified through generators this factor with a component of $\mathscr N$ (from our list).  (This component is uniquely determined since the (BS)invariants $I_{C,C'}^s$ are multi-linear, so the factors are given by pairwise distinct co-ordinate vectors.)  Then taking $i-1=\textbf{g}$, we conclude by Krull's theorem that every component of $\mathscr N$ is in our list, in other words that the component map is surjective.

  There is of course a flaw in this argument, which we call the ``elephant in the room'', because one may cannot easily ignore it, ``In spite of all temptations''. This can take two forms.

  First components can ``merge'' on linearisation. However this \textit{cannot} happen here because the above factorisation of $I_i$ is into products of linear functions in \textit{distinct} sets obtained from the (BS) invariants being multi-linear.  By the same token the components of $\mathscr N$ occur with \textit{multiplicity one}.

  Second a component may ``disappear'' on linearisation.  This can also be easily ruled out if linearisation involves only homogeneous substitutions because $\mathcal I$ is a graded ideal.  Now in fact our Weierstrass section is defined by ``excluded roots'' which corresponds to setting a co-ordinate vector equal to zero (a homogeneous substitution) or by the lines in a component tableau $\mathscr T^\mathcal C$ labelled by a $1$ which corresponds to setting a co-ordinate vector equal to $1$ (an inhomogeneous substitution).  Thus all we have to do is to deal with this latter issue which nevertheless threatens to be rather tricky; but we hope soon to resolve.

  \

  \subsubsection {A Possible Way to depose of the Elephant} \label {3.6.1}

  \

  ($A$).  An easy case of the above occurs for the parabolic defined by the composition $(2,1,1,2)$.  There are two invariant generators  $I^1_{C_2,C_3},I^2_{C_1,C_4}$. The first is already linear, so it produces no elephant!  Modulo this linear invariant, the second invariant factorizes as a product of $2\times 2$ minors, from whose zero sets the required components are obtained.  This simple picture becomes blurred by the composition $(3,2,2,3)$ because the first generator  $I^1_{C_2,C_3}$ is quadratic and the second $I^2_{C_1,C_4}$ of degree seven and does not factorise modulo the first.  Thus it has to be linearised producing a possible elephant.  Nevertheless there are again just two component tableaux so one may anticipate just two components and the hopefully the elimination of the elephant.

  \

  ($B$).  Consider the parabolic defined by the composition $(1,2,1,2)$.  There are two generators of $I$, namely $I^1_{C_1,C_3}$ which is quadratic and $I^2_{C_2,C_4}$ which is cubic.

  The linearisation of $I^1_{1,3}$, via the Weierstrass section of \cite [Thm. 5.2.6]{FJ5} requires setting
$$x^*_{1,2}=1, \quad x^*_{1,3}=0.$$

In this manner $I^1_{1,3}$  becomes the linear term $x^*_{2,4}$.

The first of the displayed identities, destroys homogeneity. However it is \textit{not} required to obtain factorisation of $x^2_{2,4}$, since $x_{1,2}$ is absent in this invariant.  Hence no elephant!

Then if we factor out by the resulting linear term $x_{2,4}$ which is homogeneous, the second invariant $I^2_{2,4}$ factors as $x_{3,4}(x_{2,5}x_{4,6}-x_{2,6}x_{4,5})$ giving two components by Krull's theorem.

 Moreover we can also obtain this factorisation by \cite [1.10]{FJ2} because the corresponding reverse tableau, obtained by placing $4$ under $2$ and pushing $3$ under $1$, has three columns of height $2$.

 This further results in two reverse tableaux obtained by either by shifting $4$ into $R_3\cap C_1$ or by shifting $6$ into $R_3\cap C_2$. Finally the red multi-sets are respectively $(4,4)$ and $(4,6)$ and can be thereby identified with the two components coming from the component tableaux, incidentally proving surjectivity of the component map in this (baby) case. This is illustrated in Sect. \ref {10}, Example $2$.

 This provides a paradigm for the general case though we cannot expect to dispense so easily with the inhomogeneous term (and hence the elephant).  Already the second example in \ref {4.4.4} is simple but tricky, and which we cannot so far resolve.

\section {Reverse Tableaux}\label {4}

\subsection {Reverse Tableaux} \label {4.1}

Recall the notion of a component tableau \cite [1.7]{FJ5}.  Here we shall construct a new family of tableaux called Reverse Tableaux, starting at the first step with $\mathscr T$.

Recall the notion of a complete sequence  (\ref {3.5}) of pairs of neighbouring columns of $\mathscr D$.  Fix one such sequence with $j^{th}$ term the neighbouring pair  $P_j:=(C_{r_j}(s_j),C_{r'_j}(s_j))$ of height $s_j$.  In our present work we will have to consider more than one such sequences but not in general all of them.

A column in a reverse tableau $\mathscr R$ is just a column in $\mathscr D$ but with entries given by $\mathscr R$.

Fix some pair $P_i:=(C,C')$ of height $s$ in this sequence.

Let $\mathscr R^i$ denote one of the several reverse tableaux
obtained by implementing (using the procedure given below) the first $i-1$ pairs in the sequence,
so in particular $\mathscr R^1=\mathscr T$.  We shall see \ref {4.3} that for $i>1$, the numerical entries of $\mathscr R^i$ may be shifted leftwards from their positions in $\mathscr T$, some entries may be repeated in red, whilst for each numerical value exactly one appears in black.


\

\textbf{Definition}.  A column of a reverse tableau is said to be of black height $s$ if its lowest black entry lies in $R_s$.

\

Thus the height of a column in a reverse tableau is its black height plus the number of red entries below its lowest black entry.

\subsubsection {The Left Boundary Column} \label {4.1.1}

View $(C,C')$ as columns in $\mathscr R^i$ and let $\mathscr I_i$ denote the set of entries which lie in the columns of $[C,C']$.  Let $C^-$
be the rightmost column in $\mathscr D$, so that
at least one entry of  $\mathscr I_i$ \textit{either black or red}
lies in the columns of $[C^-,C']$. It is called the left boundary column.  Of course it depends on $i$ but since the latter is fixed we omit this label.

\subsection {The Enabling Proposition (for the pair $(C,C')$ relative to $\mathscr R^i$)} \label {4.2}

Consider the set $S_i$ of columns of $\mathscr D$ with entries given by  $\mathscr R^i$ in $[C^-,C']$ of black height $\leq s$ and height $\geq s$.

\begin {prop} At most the leftmost column of $S_i$ has black height
$<s$ and moreover it has height $s$.

\end {prop}

The proof, due to the first named author, will be given in the sections below.  Its difficulty lies in the fact that entries change from a reverse tableau to its sequel, as pairs of neighbouring columns are implemented.

\subsection {Construction of the Reverse Tableau and Reverse $j$-strings.} \label {4.3}

\

Assume that (one of) the reverse tableaux $\mathscr R^i$ has been constructed and retain the notation of \ref {4.1}.  We assume the enabling proposition for the $\mathscr R^i$ has been proven with respect to  the $P_j:j<i$ implemented successively.  Then we implement the pair $P_i=(C,C')$ to obtain further reverse tableaux  $\mathscr R^{i+1}$, as in $(*)$ below (there being several possible choices).

\

\subsubsection {The Construction of  $\mathscr R^{i+1}$}\label {4.3.1}

\

Choose any column $C''$ of $S_i$ which is not the leftmost one.   Then by the enabling proposition for $\mathscr R^i$, its lowest entry, say $j$ is in black and occurs in $R_s$.  Recolour it in red and move the original black entry  into row $R_{s+1}$ of the first column $C'''$ having height $\geq s$ of $S_i$, strictly to left of $C''$.  \textit{In particular this black entry is moved into $[C^-,C']$, and strictly below $R_s$.}  Different choices of the columns in $S_i$ will lead to different elements of $\mathscr R^{i+1}$, but being parsimonious with notation we shall not denote them differently. (It is not even clear what would be a reasonably efficient way of denoting them.)

We call this procedure ``Substitution''.

\


Substitution will not create a gap in $C'''$, but the box in question, namely $R_{s+1}\cap C'''$ may be already occupied!  To make room for the new black entry, it will like an uninvited cuckoo-chick, starting from $C'''$, sequentially shift in unison the lower parts of columns, that is to say parts strictly below $R_s$, to the left, skipping over intermediate columns of height $<s$ to avoid creating a gap.

We call this procedure ``The Shifting of Partial Columns''.  (This is exactly analogous to ``shifting'' introduced in \cite [4.1.3]{FJ5}, for component tableaux.)

It implements a ``Skewing to the Left'', that is to say each row is moved more and more to the left as we go down columns, corresponding to the different possible values of $s$.

 Notably a column in $\mathscr R^1$ will appear in $\mathscr R^i:i>1$ skewed to the left, that is to say by joining their entries we will obtain  a line (not quite straight!) moving to the left going downwards.

\subsubsection {Extreme Shifting}\label {4.3.2}

\

This process of shifting entries below a certain row to the left may be allowed to stop, without creating gaps in columns, if a column of height $s$ is reached so in particular at the leftmost column of $S_i$ which by the Enabling Proposition has height $s$.  In this, we may or may not choose to skip over intermediate columns of height exactly $s$.  Fully skipping over is referred to as Extreme Shifting.  This was used  to obtain a useful inclusion \cite [6.3]{FJ6}, see also \ref {9.4}. Generally speaking it increases the set of excluded roots in the reverse tableau - see Example $4$ as well as Example $5$ of Section \ref {10}.  In the latter case we show in Section ref {9}, that these extra roots are ``covered'', so by Theorem \ref {9.3}, this does not change the variety defined by excluded roots.

The use of extreme shifting will be circumvented here by Proposition \ref {5.2}.

\subsubsection {Reverse $j$-strings}\label {4.3.3}

\

When $m'_{r+1}\leq 1$ the corresponding reverse tableau and of course the last one when all pairs in a complete sequence have been used will, for some numerical entry $j$, admit $m_{r+1}'$ red entries and then a last black entry, all with numerical value $j$.  These coloured entries with the same numerical value $j$ go strictly  to the left starting from the original black entry (now recoloured in red) and successively down by exactly one row, ending with the last black entry,
properties not changed by any subsequent shifting of lower parts of columns (to the left).

\

We call this a reverse $j$-string.  It consists of all the entries with numerical label $j$. Notice that as a consequence the following holds.


\begin {lemma}

\

$(i)$.  For all $j\in [1,n]$ there is at most one box with numerical entry $j$ of whatever colour in a given row of a reverse tableau.

\

$(ii)$, The lowest appearance of $j$ is in black and successively in one row above they are in red till the row in which the original black entry, that was recoloured in red,  is reached.

\

$(iii)$. If $k \in C^-$ (coloured either in black or red), then it is the unique entry of $k \in [C^-,C']$.

\

$(iv)$.  Column-wise the red entries in $\mathscr T$ lie strictly between the original black entry and the leftmost black entry.

\

$(v)$.  For each numerical value there is exactly one black entry and they form the black set $B(\mathscr R)$ of the reverse tableau $\mathscr R$.

\

$(vi)$.  By contrast for each numerical value there may be none coloured in red or several.
\end {lemma}

\begin {proof} Clear from the comments preceding the statement of the lemma.
\end {proof}

\textbf{Remark $1$}.  The colouring in the reverse $j$ - strings is irrelevant in the sense that there is a unique leftmost entry of a given numerical value which we can then declare to be in black and the remaining entries (if any) to be in red.

\

\textbf{Definition}. A tableau $\mathscr R$ is said to be a union of reverse $j$-strings, if its set of numerical entries lie in reverse $j$-strings (consisting of all the entries with numerical label $j$) as $j$ runs over distinct (not necessarily all) values in $[1,n]$.

\

In particular a reverse tableau is a union of reverse $j$ strings, some of which may be singletons.

For examples, see Section \ref {10} paying attention to the fact that some figures are of component tableaux, so do not have the properties cited above.

Adopting the colouring prescribed by Remark $1$, set of red entries in which a numerical value may appear many times.  It is called the Red Set of the (reverse) tableau and denoted $R(\mathscr R)$.  Properly speaking it is also a multi-set.

The  multiplicities in the appropriate lowest elements of a given column of $\mathscr T(\infty)$ will coincide with the appropriate multiplicities of the  Red Set for the corresponding reverse tableau $\mathscr R^{\psi(\mathcal C)}$ - to be constructed in \ref {6.1}.  In this there is a significant change of convention.  In a component tableau, a red entry with a given numerical value occurs only once, but is assigned a multiplicity, given by the number of vertical lines labelled by $\ast$ which it meets.  In a reverse tableau there may be several boxes with a given numerical entry, lying below the second end-point of the vertical lines coming from the component tableau.

A further step in our Royal Route is the following

\

\textbf{Exercise $1$}.  Compare the Red Set $(4,4)$ of the component tableau in Example $2_1$ which is a multi-set with up-going vertical lines from $4$ to $3,2$ with the corresponding reverse tableau with the same Red Set given in the top half of Figure $2$ of Section \ref {10}, in which there are \textit{two} copies of $4$ one of which goes upwards to $2$ and the second to $3$.

\

\textbf {Remark $2$.} Extreme shifting produces a unique reverse tableau having a given Red Set (defined by choices within each $S_i$ in which black entries are recoloured in red).  However this reverse tableau still depends on the choice of the complete sequence.  Thus one cannot expect to obtain \textit{naturally} a unique reverse tableau for a given Red Set, unlike the case of the component tableaux which are uniquely determined by their Red Set, via Lemma \ref {3.3}.  Examples are given through Figure $10$ of Section \ref {10}.

\

A further step in our Royal Route to tableaux is the following

\

\textbf{Exercise $2$}.  Compute the reverse tableaux associated to the component tableaux of Figure $1$, Sect. \ref {10}.  The answer ia given in the right hand column of  Figure $10$, Sect. \ref {10}; but there is some ambiguity heere as the reverse tableaux associated to a given Red Set are not uniquely determined, for example if Extreme Shifting is used and in what order pairs of neighbouring columns are implemented.


\

The rigidity of Component Tableaux implies that they are uniquely determined by their Red Set via Lemma \ref {3.3}, hence ``reliable''.  Yet the contrasting mobility of Reverse Tableaux is crucial to factorisation.  This may be understood by examining the right hand column of Figure $10$ of Section \ref {10}, in which two different reverse tableau with the same Red Set $\{7,8,8,11\}$ are obtained by changing the ordering of the implementation of the last two pairs $(C_2,C_7),(C_3,C_6)$ necessitated for factorisation.  Yet these define the same variety through excluded roots (Theorem \ref {9.3}). A further major advantage of Reverse Tableaux over Component tableaux is that it much easier to extract the excluded roots in the former (\ref {4.3.5}) and so to obtain the components they define, which only depends on their Red Set (\ref {9.3}).
Again the proofs of \ref {8.1} are much easier for Reverse Tableaux. In this sense, Reverse Tableaux  are ``user-friendly'', but their ever-changing nature might make them seem ``deceptive'' and difficult for some to
handle\footnote{As the song goes:  Sempre un amabile/ Leggiadro viso/In pianto o in riso/È menzognero.} - indeed even for their existence we need the Enabling Proposition whose proof could only be provided by the first named author \ref {4.4} by ``moving the goal posts''.

\subsubsection {A Structure Property of reverse Tableaux} \label {4.3.4}

\

\textbf{Definition}. We say that a tableau is standard with multiplicities, if entries increase strictly down rows (standard for columns) and strictly from left to right (standard for rows), with possibly more than one box with the same entry.

Notice that $\mathscr T$ is standard (with no multiplicites) and even in the stronger sense that if $i\in C_r$ and $j\in C_{r+1}$ then $i<j$.

Recall the notation of \ref {4.1}; but replace the index $i$ by $u$ in this subsection.

\begin {lemma}  A Reverse tableau is standard with multiplicities.
\end {lemma}

\begin {proof}

We prove the assertion by induction on $u$.  It holds for $u=1$ since $\mathscr R^1=\mathscr T$ is standard (with no multiplicities).

Now apply the construction of \ref {4.3.1} to $\mathscr R^u$. This involves substitution, which recolours in red a lowest black entry $j$ in some column $C'$ of black height $s$ and places $j$ in the first column $C$ of height $\geq s$ strictly to the left of $C'$ in row $R_{s+1}$.

Obviously recolouring does not upset the standard property.

Again in the above $j$ is strictly greater than the entry $R_s\cap C$, by $\mathscr R^u$ being standard with respect to rows and so strictly greater than all the entries of $R^s\cap C$, by $\mathscr R^u$ being standard with respect to columns.

Now observe that the property of being standard with multiplicities is also preserved by the Shifting of Partial Columns.  This derives from the fact that given an entry $i$ in a standard tableau, moving to an entry obtained by going horizontally to the right and then downwards to an entry $j$ implies that $i<j$.  Combining this with the fact that the lower parts of columns are moved in unison to the left with entries staying in the same row, implies the required result.

\end {proof}

\textbf{Example $2$}.  Consider the composition $(1,2,1,2)$.  In the first step $4$ is moved under $3$ in $C_2$ and is recoloured in red in $C_3$. Obviously this satisfies the standard property with the two entries of $4$ being in different rows and columns.

Now in the second case we have two possible new tableaux - See Figure $2$ of \ref {10}.

In the first possibility, $4$ is moved below $3$ in $C_1$, so the standard property is preserved; but with two copies of $4$.

In the second possibility $6$ is moved under $4$ in $C_2$ strictly to the left but since it is is strictly greater than $4$ the previously moved entry, the standard property is again preserved.

\subsubsection {Excluded Roots in the Reverse Tableaux}\label {4.3.5}

The (defined) set of excluded roots $X(\mathscr R^\mathcal E)$ for a reverse tableau is extremely simple indicating the advantage of the latter, but not the most important one.

First recall that a reverse tableau has entries which increase strictly along rows and down columns, though if one counts red entries some numerical entries may be repeated (\ref {4.1}).

\

\textbf{Definition}. Define $x_{i,j}:i<j$ to be excluded, that is in $X(\mathscr R^\mathcal E)$, if $i$ appears above $j$ in the same column or in a column strictly to the right of the column containing $j$, and only if $i$ is the rightmost choice of $i$ in its unique reverse $i$  (By convention roots in the Levi are not counted, as in \cite [4.1.4, definition]{FJ5}).

\

 Suppose $i$ appears to the right of $j$, that is the line $\ell_{j,i}$ is right going.

  If $i$ also occurs in a row (not necessarily strictly) below that containing $j$,  then $i>j$, by the standard property. Thus the condition that $i<j$ implies that $i$ always lies in a row strictly above that containing $j$.  Consequently we have the

  \begin {lemma} If a right going line $\ell_{j,i}$ is downward going (for example horizontal) then $x_{i,j}$ is not an excluded root.
  \end {lemma}

\textbf{Remarks}.  Thus a horizontal line does not give rise to an excluded root.  Again in the above we need only take $j$ to be coloured in black, since any red entry always lies strictly to the right of the (unique) black entry having the same numerical value.  On the other hand we sometimes must allow red to be a colour of the first entry $i$.

\

\subsubsection {Excluded Roots in the Component Tableaux}\label {4.3.6}

\

In a component tableau $\mathscr T^{\mathcal C}$ we do not obtain its set $X(\mathscr T^{\mathcal C})$ by reading off from the columns of the tableau. Rather for each red line $\ell^\mathcal C_{i,j}$ with $i\in \mathscr T^\mathcal C(\infty),j\in J^\mathcal C_{r+1}$ we create \cite [4.1.3]{FJ5} a new tableau $\mathscr T^{\mathcal C}_{i,j}$ (the superscript $\mathcal C$  does not appear in loc.cit.) by placing $j$ directly under $i$ as it appears in $\mathscr T$ and shifting lower parts of columns to the left to make room for $j$. Again although couched in terms of the Weyl group the set $X(\mathscr T^{\mathcal C}_{i,j})$ of excluded roots of $\mathscr T^{\mathcal C}_{i,j}$ is exactly as defined for the reverse tableau as above.  Finally $X(\mathscr T^\mathcal C)$ is obtained by taking the union of the excluded roots coming from each red line \cite [4.1.5, Definition]{FJ5}.

\begin {lemma}  $\ell^\mathcal C_{i,j} \in X(\mathscr T^{\mathcal C}_{i,j})$.

%
\end {lemma}

\begin {proof} This obtains since $j$ is placed below $i$ to obtain $\mathscr T^{\mathcal C}_{i,j}$ and through the definition of $X(\mathscr T^{\mathcal C}_{i,j})$. In slightly different language it is \cite [Lemma 4.5]{FJ5}.

\end {proof}

\textbf{Remark $1$.} Again we shall define (\ref {6.1}) a reverse tableau $\mathscr R^{\psi(\mathcal C )}$ with the same Red Set as $\mathscr T^\mathcal C$.  Now let $\ell_{i,j}^{\psi(\mathcal C)}$ be $\ell^\mathcal C_{i,j}$ but viewed as element of $\mathscr R^{\psi(\mathcal C )}$.  We showed in \cite [4.1.3]{FJ6} that $\ell_{i,j}^{\psi(\mathcal C)}\in \mathscr R^{\psi(\mathcal C )}$. At first sight this seemed rather tricky and we used the appearance of the amalgamated column (\ref {6.1.1}) in the reverse tableau.  A direct proof is given in ref {8.1.1} using the rightmost appearance of $i$ in the reverse tableau, which makes it ``easier'' for $\ell_{i,j}^{\psi(\mathcal C)}$ to be an excluded root for the latter.

\

\textbf{Remark $2$}.  The amalgamated column could also have been used in the proof of $(i)$, which simultaneously treats all the excluded roots $\alpha_{i,\cdot}$, for all the lines $\ell^\mathcal C_{i,\cdot}$, labelled by a $\ast$.

\

\subsection {The Trapezium} \label {4.4}

\subsubsection {The Left Boundary\label {4.4.1}}

Fix a neighbouring column pair $C.C'$ and let $s$ denote their mutual height.  let  $i\in [1,n+1]$, be uniquely determined by the condition that  $C=C_{r_i}(s),C'=C_{r'_i}(s)$, with $s_i=s$.   We keep this notation till the end of Section \ref {5}.

 The rectangle $R^s_{C,C'}$ is modified successively by implementing the first $j-1:j\in [1,n+1]$ pairs of neighbouring columns in the complete sequence of \ref {3.5}, to obtain the $j^{th}$  trapezium $\mathcal T^s_{C,C'}(j)$, based on the $i^{th}$ neighbouring pair $C,C'$ of height $s$.

 This means that we start with $\mathcal T^s_{C,C'}(1)=R^s_{C,C'}$, and successively implement $P_1, P_2,\ldots,P_{j-1}$, to obtain $\mathcal T^s_{C,C'}(j):j \in [1,n+1]$.  Notice that here if $j\leq i$ the pair $(C,C')$ is \textit{not being} implemented.

 \

 \textbf{N.B.} Recall \ref {4.1}. The resulting Trapezium is a sub-tableau of (one of the) reverse tableaux $\mathscr R^i$  and correspondingly there may be several obtained at a given step but we do not denote them differently so as not to introduce cumbersome notation. Indeed we already have two sets of indices $(C,C')$ and $i$, decorating the Trapezium.


Recall that a given numerical entry of some $\mathscr R^j$ may be either in black or in red.

\

For all $t \in [1,s]$, let $b^j_t$  be the rightmost box in $R^s$, \textit{ containing in $\mathscr R^j$}, the numerical entry, either in black or red, of $R_t\cap C$.


\

\textbf{Remark.}   Our construction \ref {4.3.1} recolours the rightmost entry having some numerical value in red and places a black entry with that numerical value strictly to its left. Thus $b^j_t$ being the rightmost entry, belongs to row $R_t$, but by shifting may appear in a column to the left of $C$.


\

\textbf{Definition}. The left boundary $B^j$ of $\mathcal T^s_{C,C'}(j)$, is just the set of the $s$ boxes $\{b^j_t\in R_t\}_{t=1}^s$.

\subsubsection {The Right Boundary\label {4.4.2}}

For all $t \in [1,s],j \in [1,i]$, let $b'^j_t$  be the leftmost box in $R^s$ \textit{containing in  $\mathscr R^j$}, the numerical entry of $R_t\cap C'$.

\

\textbf{Definition}. The right boundary $B'^j$ of $\mathcal T^s_{C,C'}(j)$ is just the set of $s$ boxes $\{b'^j_t\in R_t\}_{t=1}^s$.

\

However unlike the case of the left boundary it is not obvious that $b'^j_t \in R_t$. We prove this by induction on $j\leq i$ and that this entry is always black.  This is obvious for $j=1$, since the trapezium is just the rectangle. However it fails when $j>i$.



\



\subsubsection {Definition of the Trapezium\label {4.4.3}}

Again assume $j\leq i$.

The $i^{th}$ trapezium (resp. left $j^{th}$ trapezium) $\mathcal T^s_{C,C'}(j)$ (resp. $^l\mathcal T^s_{C,C'}(j)$) is defined to be $ R^s\cap [B^i,B'^i]$ (resp.  $R^s\cap ]B^j,B'^j]$).


\


By construction the numerical entries of $B^j$ (resp. $B'^j$) are the same of those of $C$ (resp. $C'$).

Those in $B^j$  are the \textit{unique numerical entries in the trapezium}.  This is also true for $B'^j$ as shown in Lemma \ref {4.4.4} below.


  Joining the entries of $B^j,B'^j$ gives (not quite straight!) downward, left going lines.  Thus the term trapezium, also bounded by the horizontal lines above (resp. below) by $R_1$ (resp. $R_s$).


\

\textbf{Example $2_2$.} Consider the composition $(3,1,3,1,3)$, with $C=C_3,C'=C_5$ and implement the complete sequence  $(C_2,C_4);(C_1,C_3);(C_3,C_5)$ in that order. The left rectangle $^lR^3_{C_3,C_5}$ has $4$ black entries, being exactly the degree of
$I^3_{C_3,C_5}$.  After implementing the pair $(C_2,C_4)$, one finds that $B$ is just $C$ skewed to the left whilst $B'=C'$.   After implementing the pair $(C_1,C_3)$,  one finds that $B$ is unchanged except that $7$ is recoloured in red and a black $7$ appears in $R_4\cap C_1$.  It is exactly at this stage that the number of black entries in the left trapezium is decreased by one. Again $B'$ is unchanged.

We do not illustrate this example, but the reader should carry out the computation following the instructions.

\subsubsection {A Hidden Rule\label {4.4.4}}

 Implementing the neighbouring pair $C=C_{r_i},C'=C_{r'_i}$ on the Trapezium $\mathcal T^s_{C_{r_{j}},C_{r'_{j}}}(j)$
 obtained from the $j^{th}$ step, is to be applied \textit{only} to the entries of this trapezium. 
 This condition, as natural and innocuous as it might look, actually automatically implements  the condition
 that a black entry which is not subject to substitution by this rule is exactly one for which the \textit{original} black entry was not surrounded by the pair $(C,C')$.  A similar restriction occurs for the component tableau and furthermore there are correspondingly less components  resulting from factorisation.

 This is illustrated in the following example.

 \
%
%
%

  \textbf{Example $2_3$.} Consider just the composition $(1,2,2,1)$.    Implement the pair $(C_1,C_4)$. This recolours the $6$ in red and places a black $6$ under $4$, and then $5$ under $2$ and $3$ under $1$. As above we cannot further recolour the new black $6$ in red by our hidden rule.  Again there is \textit{no} component tableau with red multi-set $(6,6)$.

 Finally the evaluation that this implements on $I^2_{2,3}$ puts a zero at $x_{2,5}$ which does \textit{not} eliminate components because it is a homogeneous substitution whilst this invariant factors as $x_{2,4}x_{3,5}$.
 On the other hand this evaluation of $I^2_{2,3}$ does not depend on any particular value given to $x_{3,5}$ which could be set equal to zero leaving only the one component for $\mathscr N$.

 Thus one suspects that $\mathscr N$ is irreducible in this case.  It is one of the simplest examples for which one should try to get rid of the elephant.

 \begin {lemma}  Retain the notation and hypotheses of \ref {4.4.2}, in particular $j\leq i$. Then $b'^j_t \in R_t$. It is the unique entry of $j^{th}$ trapezium) $\mathcal T^s_{C,C'}(j)$ with numerical value that of $R_t\cap C'$ and is in black.  Thus the right boundary $B'^j$ forms a left going downward line of black entries, coinciding row for row with those of $C'$, which appears thereby skewed to the left.
 \end {lemma}

\begin {proof}  Suppose we have implemented a neighbouring pair $C_1,C_1'$ of height $t$ on $\mathcal T^s_{C,C'}(j-1)$.  If $t>s$, this will not change the entries in this trapezium.  Since we are not implementing the pair $C,C'$ by hypothesis, we only remain with $t<s$.

This can only recolour a black entry in $R_t\cap C'$ in red (by substitution) if the column containing this element has black height $t$.  This must mean that the original black element in $R_{t+1}\cap C'$ has been shifted horizontally to the left and replaced  by a red element in a right adjacent column (it being of no consequence if replaced by a new black entry, which could also happen).

The latter operation would require a neighbouring pair $(C_{(2)},C'_{(2)})$ surrounding $C'$ of height $t$.  However in order to further recolour the black entry in $R_t\cap C'$ in red, we would need to implement a second column pair $(C_{(3)},C'_{(3)})$ (\textit{different} to $(C_{(2)},C'_{(2)})$) also of height $t$, surrounding $C'$. This is impossible as $C'_{(3)}$ would then need to be to the left of $C_{(2)}$ and so strictly to the left of $C'$.
\end {proof}

\textbf{Remark}.  We see here the importance of not implementing a given neighbouring pair twice.  However let us again emphasize that this restriction is only important for the right boundaries $B'^j$.

\

\textbf{Example $2_4$}.  Consider the composition $(2,1,3,1,3,2,1)$ with $C=C_3,C'=C_5$. Implement the pair $(C_1,C_6)$ having height $2$. This places a red $12$ below the black $11$. Now implement the pair $(C_{(2)}:=C_4,C'_{(2)}:=C_8)$ of height $1$. Then a black $13$ pushes $C_5^{>1}$ to the left and places a red $12$ under the black $8$, with $C'$ now of black height $1$. Here $B'$ is skewed to the left, with its highest entry the black $8$ which is now in a column of black height $1$. Yet we cannot recolour this black $8$ in red by substitution and indeed the only available pair of height $1$ is $(C_{(3)}:=C_2,C'_{(3)}=C_4)$ whilst $C_4$ lies strictly to the left of $8$, itself lying in $C'=C_5$.

\subsubsection {Non-overlapping of trapezia of the same height}\label {4.5.5}

In analogy with the essential non-overlapping of pairs $P_{i_1}:=(C_1,C'_1);P_i:=(C,C')$ of neighbouring columns of the same height $s$, we have the following

\begin {lemma} Distinct trapezia $\mathcal T^s_{C_1,C'_1}(j),\mathcal T^s_{C,C'}(j):j\leq i_1$, appearing in that order, that is with the second on the right, can at most overlap by sharing a common boundary.
\end {lemma}

\begin {proof} The rectangle $R^s_{C_1,C'_1}$ must lie strictly to the left of $R^s_{C,C'}$ with either in the first case $C$ strictly to the right of $C'_1$ or in the second case equal to it.  Now by choice of $j$, the right boundary of the left hand trapezium $B_1'^j$  has the same numerical entries as $C_1'$, whilst by Remark \ref {4.4.1}, the left boundary of the right hand trapezium $B$ has the same numerical entries as  $C$.

In the first case (because a tableau is standard, that is to say with entries stricly increasing down columns and from left to right in rows) the entries of $C_1'$ are strictly less than those of $C$.  On the other hand a reverse tableau is standard with multiplicities \ref {4.3.4}, in particular  entries must increase from left to right in rows.

Thus in this case $B'_1$ must lie row by row strictly to the left of $B$.

  In the second case we claim that the entries of $B'_1$ and of $B$ coincide  in numerical value and in colour.  Indeed the numerical entries of $B_1'$ (resp $B$) coincide row by row with the numerical value in the common column $C_1'=C$ and so coincide in colour by Lemma \ref {4.3}(i).

  In brief one may identify $B'_1$ and $B$.
  \end {proof}

%
%
%
%
%
%

  \

  \textbf{Remark}.  Of course in the first case $(C'_1,C)$ is also a pair of columns of height $s$, though not necessarily neighbouring.  Thus why consider Case $1$ at all!  This arises because we want to choose an order in implementing neighbouring columns, which regulates implementing the pair $(C'_1,C)$ to say a later stage.

\subsubsection {Non-acquisition of Red Entries}\label {4.4.6}

\begin {lemma} The left trapezium $^l\mathcal T^s_{C,C'}(j)$ can only acquire a red entry in $R_s$ by implementing the pair $C,C'$.
\end {lemma}

\begin {proof} First consider substitution.  In this the black entry $b$ to be recoloured must lie in $R_s$. Its recolouring can only come from a neighbouring trapezium $\mathcal T^s_{C_1,C'_1}(j):j\leq i_1$ on the left of $\mathcal T^s_{C,C'}$, if they share a common boundary $B$ being also of course the left boundary of $\mathcal T^s_{C,C'}(j)$.

 This however results in a black entry of $b\in B\cap R_s$ being recoloured in red, but of course the resulting red element remains in $B$ and so  is not in the left trapezium  $^l\mathcal T^s_{C,C'}(j)$.

Second, consider translation.  In this case one must start with a red entry in $R_s$ strictly to the right of $^l\mathcal T^s_{C,C'}(j):j\leq i$.  However it is then blocked by the right boundary of this trapezium which can only have black entries.

\end {proof}

\subsubsection {Non-interference of Columns of Height $>s$.\label {4.4.7}}

Observe that the  implementation of a pair $P_j:j<i$ of height $s'\geq s$ does not modify the numerical entries of the trapezium $\mathcal T^s_{C,C'}(i) $, by substitution because a black element in row $R_{s'}:s'\geq s$ is recoloured in red, whilst a black element in $R_{s'+1}\cap [C_{r_i}^-,C_{r'_i}]$ is adjoined. \textit{However} if a column to the left of the pair $C,C'$ of height $s$ is implemented, then an entry of $B^j$ may be recoloured in red.

Again the trapezium $\mathcal T^s_{C,C'}(i) $ is not modified by the implementation of a pair  $P_j:j<i$ of height $s'\geq s$ through shifting in the sense of \ref {4.3}, since only rows of strictly below $R_{s'}$, so strictly below $R_s$, are shifted.

Finally observe that, as an essential consequence of ``moving the goal posts'', the entries of a trapezium are not modified by shifting, and the entries of the boundaries also remain though may move  leftwards.  Finally internal indices stay internal though may be duplicated in red by forming a reverse string.

  This behaviour is exemplified in Example $2$, with a black $8$ entering $R_2\cap C_3$ with $B$ obtained by skewing $C$ to the left.

  \

  $(*)$.  To summarize, in implementing the pairs $P_j:j<i$, the entries in the left rectangle $^lR^s_{C,C'}$ still appear in the left trapezium $^l\mathcal T^s_{C,C'}(i) $ - see also the last part of Proposition \ref {5.1.2}.

  \

  \textit{This is the key combinatorial advantage of trapezia.}

 \subsubsection {Proof of the first part of the Enabling Proposition.\label {4.4.8}}

  by \ref {4.4.1}, the left boundary $B$ has its lowest entry in $R_s$, and this may be in black or red.  Moreover the column $C^-$ of $\mathscr R^i$ in which the lowest entry of $B$ occurs has height $s$ and black height $\leq s$, the latter being possibly a strict inequality.  Furthermore $C^-$  is the rightmost column of $\mathscr R^i$ such that $R_s\cap C^-$ has the same numerical entry as $R_s \cap C$ in $\mathscr T$.

\

$(*)$. \textit{The observations of this last paragraph, proves the first part of the enabling proposition \ref {4.2}.}

%

\subsection {Proof of Enabling Proposition}{\label {4.5}

\

By \ref {4.4.8}, we have only to show the last part.  By Lemma \ref {4.4.6}, there are no red entries in $R_s\cap \ ^l\mathcal T^s_{C,C'}(i) $.  This means that any column in $S_i$ of height $\geq s$ and black height $<s$ must have a red entry in $R_s$ and so can only be the leftmost column of $S_i$.

\section {Vanishing of the Benlolo-Sanderson Invariants}\label {5}

\subsection {Goal}\label {5.1}

Recall that a complete sequence $P_i:=\{(C_{r_i}(s_i),C_{r_i}'(s_i))\}_{i=1}^{\textbf{g}}$ of pairs of neighbouring columns of height $s_i$, defines a linearly ordered sequence of Benlolo-Sanderson invariants $I_i:=I^{s_i}_{C_{r_i},C'_{r_i}}; i\in [1,\textbf{g}]$.

Set $\mathscr R^1=\mathscr T$ and inductively define $\mathscr R^{i+1}$ to be  any reverse tableaux constructed from any $\mathscr R^{i}$ by implementing the pair $P_i$
using the enabling proposition, as detailed in \ref {4.3}${(*})$.  \textit{Again} at each step there may be several different choices and therefore several different ``final'' reverse tableaux which we are not denoting differently (indeed this would require cumbersome notation).


Let $\mathfrak u_i$ be the span of the roots in $\mathfrak m$ which do not belong to some $X(\mathscr R^i)$ and denote $\mathfrak u_{\textbf{g}+1}$ simply by $\mathfrak u$.  We eventually wish to show that $\overline{B.\mathfrak u}$ is a component of $\mathscr N$.  Again it depends on choices and so on the resulting ``final'' reverse tableaux. Consequently we may thereby obtain several, in general different, components of $\mathscr N$.

The basic step is the

\begin {thm} The Benlolo-Sanderson invariant generators  vanish on $\mathfrak u$ and so on $\overline {B.\mathfrak u}$.
\end {thm}

This will be established in the sections below.

\subsubsection  {A Consequence of \ref {4.3}. }\label {5.1.1}

\

Suppose that a previous pair $P_j=:(C_j,C_j'):j<i$ had been implemented with $s_j<s$.  Through the enabling proposition applied at this step, this recolours in red, a black entry $k$ which occurs in $R_{s_j}$ as the lowest black entry of some column in $\mathscr R^j$ and places a new black entry $k$ in $R_{s_j+1}$ strictly to the left of the black entry which had been recoloured in red.

Recall that we are setting $C=C_i$ and $C'=C_i$, here and in \ref {5.1.2}.

\begin {lemma}  This new black entry lies in $^l\mathcal T^s_{C,C'}(i)$.
\end {lemma}

\begin {proof}

This follows inductively from the penultimate sentence of the first paragraph of \ref {4.3.1} (which covers substitution) and \ref {4.4.7}$(*)$ (which covers shifting).  It is here, in particular, we see the advantage of moving the goal posts (that is to say the boundaries of the trapezia as they evolve).

\end {proof}


\subsubsection {A Key Proposition\label {5.1.2}}

\

\begin {prop}  There are no red entries in $R_s\cap \ ^l\mathcal T^s_{C,C'} $ if and if only if the number of black entries in $^l\mathcal T^s_{C,C'}$ is the degree $d$ of $I^s_{C,C'}$.  In addition these have the same set of numerical entries as the left rectangle $^lR^s_{C,C'}$ baring the possible red entry in $R_s\cap \ ^l\mathcal T^s_{C,C'} $.
\end {prop}

\begin {proof}  Take $j \in [1,i+1]$ and let $N_{\mathscr R^j}(R^s\cap ]B^j,B'^j])$ be the number of black boxes in $R^s\cap ]B^j,B'^j]$ with respect to a reverse tableau $\mathscr R^j$ (with a corresponding meaning to this expression for individual rows, that is when $R_t:t \in [1,s]$ replaces $R^s$). By Lemma \ref {2.7}, $N_{\mathscr R^j}(R^s\cap ]B^j,B'^j]) =d$ if $j=1$.

Now suppose $N_{ \mathscr R^j}(R^s\cap]B^j,B'^j])=d$ for a reverse tableau $\mathscr R^j$, for some $j\in [1,i-1]$.  Suppose there is a box with black entry $k$  in $R_i\cap ]B^j,B'^j]$. Substitution would recolour this entry in red and introduce a black entry with numerical entry $k$ in $R_{i+1}$.  Consequently
$$N_{\mathscr R^{j+1}}(R_i\cap ]B^{j+1},B'^{j+1}])=N_{\mathscr R^{j}}(R_i\cap ]B^{j},B'^{j}])-1,$$
whilst by Lemma \ref {5.1.1} one also has
$$N_{\mathscr R^{j+1}}(R_{i+1}\cap ]B^{j+1},B'^{j+1}])=N_{\mathscr R^{j}}(R_{i+1}\cap ]B^{j},B'^{j}])+1.$$

Suppose $t>i=s$, then by taking the sum of these two expressions it follows that \newline $N_{\mathscr R^{j+1}}(R^t\cap ]B^{j+1},B'^{j+1}])=d$, whilst $N_{\mathscr R^{j+1}}(R^s\cap ]B^{j+1},B'^{j+1}])=d-1$.

The last part is an immediate consequence of the above computation.

\end {proof}

\subsection {Non-Vanishing} \label {5.2}

Recall the above notation.

Our goal is to prove the following
\begin {prop}   For all $j\leq i+1$, there are no red entries in $R_s\cap \ ^l\mathcal T^s_{C,C'} $ if and if only if $I^s_{C,C'}|_{\mathfrak u_j}\neq 0$.
\end {prop}

This would seem an almost obvious conclusion of \ref {5.1.2} and equally that Theorem \ref {5.1} would result.  These assertions are proved below.

Consider some stage of the corresponding left trapezium $^l\mathscr T^s_{C.C'}$ as it evolves from the left rectangle and let $\mathfrak u$ denote the span of the roots in $\mathfrak m$ which it does \textit{not} exclude at that stage.  Let $d$ denote the Benlolo-Sanderson invariant $I^s_{C,C'}$.

In view of \ref {5.1.2} we may reformulate the above proposition as
\begin {prop} The number of black entries in $R_s\cap \ ^l\mathscr T^s_{C,C'} =d \Leftrightarrow I^s_{C,C'}|_{\mathfrak u}\neq 0$.
\end {prop}

\subsection {Outline of proof of $\Rightarrow$}\label {5.3}

\

We can assume that the pair $(C,C')$ has not yet been implemented, otherwise the number of black entries in the left trapezium $<d$.

The proof of $\Rightarrow$ will only involve the entries of the left trapezium $ ^l\mathscr T^s_{C,C'}$.  Recall that its left and right boundaries of the trapezium both admit exactly one entry in each row $R_t:t\leq s$.  Those in left boundary may be in red  (\ref {4.4.1}), whilst those in the right boundary are all in black (Lemma \ref {4.4.4}).

\

\textbf{Definition.} A set $\mathscr L$ of right-going lines $L$ joining the black entries of $\mathscr T^s_{C,C'}$ (including the red entries of the left boundary), is said to be disjoint if for distinct lines  $\ell_{i,j},\ell_{i',j'}$ one has $i\neq i',j\neq j'$.

\

For each $L \in \mathscr L$, let $L^-$ (resp. $L^+$) be their left  (resp right) hand entries.  In this $L$ defines a bijection $\varphi_L$ of $L^-$ onto $L^+$.

\

Let $\mathscr L$ be such a set.   Then $\mathscr L$ is said to be complete if $L^-$ (and hence $L^+$) has cardinality $d$.

\subsection {The Role of Horizontal Lines}\label {5.4}

\

The advantage of horizontal lines is that they are never excluded in the reverse tsbleau  Lemma \ref {4.3.5}. Suppose we have a choice of  $\hat{L}\in \mathscr L$ in which all lines are horizontal.  We conclude that the monomial
$\prod_{i\in \hat{L}^-}\ell_{i,\varphi(i)},$ does not vanish on $\mathfrak u$.

We still have to prove that such a choice exists!  This is by induction on the number of steps  to go from the left rectangle to the left trapezium.

In the left rectangle it is certainly true.  Then each step involves substitution and translation.  The latter does not change the lines being horizonatal.

For substitution, we can assume that $(C,C')$ has not being implemented.  Then there is a black entry with numerical value $k$ in some row $R_t:t<s$ between two black entries $i,j$ in that row.  At worst $i$ is in the left boundary and $j$ is missing.

Then this entry is recoloured in red and a black entry inserted into $R_{t+1}$.  It will push an entry $i'$ to its left (in the same row) possibly a left boundary entry.  This gives a lines $\ell_{i',k}, \ell_{k,j'}$ whilst $\ell_{i',k'}$ is eliminated.  At worst $j'$ is missing.  This construction shows that the number of horizontal lines is unchanged and more significantly so are their beginning (resp. end) points of which there are $d$ in number.

 \subsection {The Interpretation of the BS Invariant $I^s_{C,C'}$ in the Trapezium} \label {1.3}

 \


A BS invariant $I^s_{C,C'}$ obtains when there are two Levi blocks of the same size, say $s$. Its description reduces to the case when these are the first and last Levi blocks and then it is the lowest degree coefficient in $c$  of the $n-s \times n-s$ minor in the bottom left hand corner restricted to $c\Id + \mathfrak m$.  When there are further Levi blocks of the same size strictly between these first and last blocks, $I^s$ factorizes to similar invariants \cite [1.10]{FJ2}.

Otherwise it is irreducible \cite [5.3]{FJ1}.  In general a BS invariant $I^s_{C,C'}$ is given by a pair neighbouring columns $C,C'$  of height $s$ in $\mathscr D$.

%
 We conclude that part of $I^s_{C,C'}$ in which co-ordinate functions with co-ordinates outside the entries of the trapezium are set equal to zero is just the determinant of linear transformations $\varphi_L:L \in \mathscr L$.

 By \ref {1.2}, the monomial coming from the horizontal lines exists and does not vanish on $\mathfrak u$, hence $\Rightarrow$.

 This does not prove $\Leftarrow$ because we have ignored part of $I^s_{C,C'}$ coming from outside the trapezium.  It is established below by a method derived from \cite [Sect. 4.3]{FJ5} (a brilliant innovation of the first named author, which avoided a seemingly long and arduous induction process).  As we shall see here \textit{further} simplification is made possible by the use of reverse tableaux.

 \

\textbf{Example.}  Return to example $2_3$.  The (BS) invariant is of degree $4$ being a product of $x_{5,8}x_{8,9}$ and $x_{6,10}x_{7,11}$.  Only the left hand factors may become excluded on implementing the pair $(C_2,C_4)$.  Take $m=5,\varphi(m)=8$ and $n=6,\varphi(n)=10$.  The said implementation puts (the black) $8$ under $5$ making $x_{5,8}$ an excluded root, as the reader was warned! Yet in this $\varphi(m)$ having pushed (on the same row) $n$ to the left, so $x_{6,8}$ is not an excluded root, nor is $x_{5,10}$.

\subsubsection {Stability of the Set of Excluded Roots.} \label {5.5.1}

Recall the definition of $\mathfrak u_i$, given in \ref {5.1}.

\begin {lemma} One has $\mathfrak u_i \supset \mathfrak u$ for all $i \in [1,n]$.
\end {lemma}

\begin {proof}  it is enough to observe that in the notation of \ref {5.1} one has $X(\mathscr R^{i}) \subset X(\mathscr R^{i+1})$.

Indeed take (one of the) component tableaux designated by $\mathscr R^i$.  Then to form (one of the) component tableaux designated by $\mathscr R^{i+1}$, we take a black entry in $j$ recolour it in red as $j_r$ and place a new black entry $j_b$ one row down and to the left. Then for any other numerical index $i$ (in any colour) $x_{i,j_b}$ is an excluded root if $x_{i,j}$ was.  On the other hand $j_r$ takes the place of $j$, so trivially $x_{i,j_r}$ is an excluded root if $x_{i,j}$ was an excluded root.  (These last excluded roots may be non-trivial ``supplementary'' excluded roots, that is by allowing their second entry to be in red.)

\end {proof}

\subsubsection {Proving ``Only if'' of Proposition \ref {5.2}} \label {5.2.5}

\

Here we proceed as in \cite [4.3.3]{FJ5}.

We need to show that the total number $N_u$ of lines forming the disjoint composite lines goes strictly down from $u=t$ to $u=t+1$ (and indeed by exactly one).

Let $c_i$ (resp $\hat{c}_i$) be the height of column $C_i$ in the trapezia $\mathcal T^s_{C,C'}(u)$ for $u=t$ and $u=t+1$.

We claim that, as in equations $(1),(2)$ of \cite[4.3.3]{FJ5}, one obtains
\[
N_t = \sum \min(s,c_i) - s, \quad N_{t+1} = \sum \min(s,\hat{c}_i) - s,
\]
where the sum is over the columns of the trapezia $\mathcal{T}^s_{C,C'}(u)$ for $u=t,t+1$.

Indeed it is obvious that disjoint composite lines can pass through $C_i$ at most $\height C_i$ number of times. On the other hand since we have exactly $s$ such lines they can pass through $C_i$ at most $s$ times if they go strictly from left to right, that is to say in common parlance that they do not ``zig-zag''.  This behaviour, for the component tableaux, was excluded in \cite [4.3.12]{FJ5}.  For reverse tableaux it is even easier.

Indeed suppose that a composite line passing though boxes, with increasing entries, were to zig-zag, then for some pair of successive entries $i,j$, the box containing $j$ would lie to the left (not necessarily strictly) of that containing $i$.

Consequently  the box containing $i$ would lie to the right (not necessarily strictly) of that containing $j$.
However since $i\leq j$, it follows from the fact that a reverse tableau is standard with multiplicities \ref {4.3.4}} that $i$ must lie in a row strictly above $j$ and then that $i<j$.  Then by \ref {4.3.5}, the root vector $x_{i,j}$ would be excluded.

Now in passing from $u=t$ to $u=t+1$ there is a column $C_i$ of the trapezium to which a black entry is adjoined to $R_{s+1}$ and a second column $C_j$ in which a black entry in $R_s$ is recoloured in red, and \textit{no} other changes.  Hence $N_t-N_{t+1}=1$ as required to prove that ``if'' implies ''only if'', and hence Prop. \ref {5.2}.

\

\textbf {Remark 1.} There is an added simplicity in using reverse tableaux since we do not need the amalgamated column of \cite [4.3.6]{FJ5}.  This may seem strange at first, but is a natural consequence of the fact that the amalgamated column appears automatically in the reverse tableaux construction. 

Again the calculations in \cite [4.3.7]{FJ5} simplify, since only one change occurs in passing from $u=t$ to $u=t+1$.

\

\textbf {Remark 2.} Notice that in the latter part of this calculation we do not need to  discount the use of red entries in the trapezium $\mathcal T^s_{C,C'}(u):u=t,t+1$, in constructing composite lines.

\subsection {Vanishing} \label {5.6}

We can now obtain a proof of Theorem \ref {5.1}.

Indeed for any Benlolo-Sanderson invariant $I^s_{C,C'}$ implementing the pair of neighbouring columns $C,C'$ after the $i^{th}$ step puts a red element into $R_s\cap \ ^l\mathcal T^s_{C,C'}(i+1) $.  Hence by ``only if'' in Prop. \ref {5.2}, we obtain  $I^s_{C,C'}|_{\mathfrak u_{i+1}}=0$.  Then by Lemma \ref {5.5.1} we obtain
$I^s_{C,C'}|_{\mathfrak u}=0$, so proving Theorem \ref {5.1}.  This is much cleaner than the proof of the corresponding result for component tableaux given in \cite [4.3]{FJ5}, and hopefully the reader is now too befuddled to disagree.

 \section {Constructing Sufficient Reverse Tableaux.}\label {6}

 Given a component tableau $\mathscr T^\mathcal C$ we construct a reverse tableau with the same Red Set which we shall denote by $\mathscr R^{\psi(\mathcal C)}$.  It is not necessarily unique but we show that it defines the same variety via excluded roots and covering, as that defined by $\mathscr T^\mathcal C$. This is proved using induction on the (natural) total order of the Red Set (\ref {3.4}), which leads to a (rather specific) complete sequence of pairs of neighbouring columns for each Red Set.

 In a short note \cite {FJ7} we show that no other Red Sets are possible for a reverse tableau.

 This indicates, but does not prove, that the component map must be surjective.  Rather the latter is to be shown by our factorisation procedure \ref {3.6}.

 Originally we had constructed $R^{\psi(\mathcal C)}$ by ``redressing'' the heights in the set of left hand members of the pairs of neighbouring columns in our chosen complete sequence, so that their heights decrease (monotonically) from left to right. Then the implementation of pairs becomes rather easier does not need ``Enabling'' and also does not involve the shifting of partial columns (which has in effect) is already implemented by redressing.

  Yet this technology is less powerful than ``Enabling'' when it comes to factorisation which requires complete sequences not defined by the order relation \ref {3.4} on the Red Set to be implemented.

  This shown can be seen from Example $14$ in Section \ref {10}.  Inspection of the  right hand column of Figure $14$ show that not all reverse tableaux with the same Red Set appear for a fixed complete sequence, whilst some appear for the same red set can be slightly different.  Thus both sequences must be used. To compare factors we must show (as we shall) that different reverse tableau with the same Red Set define the same variety through excluded roots.

 \subsection {The Construction of $R^{\psi(\mathcal C)}$ by Enabling}\label {6.1}

\begin {prop}  There exists a reverse tableau $\mathscr R^{\psi(\mathcal C)}$ with red multi-set $R(\mathscr T^\mathcal C)$.
\end {prop}

\begin {proof}  The proof is by induction on the order relation (\ref {3.4}) possessed by $R(\mathscr T^\mathcal C)$.

Assume that for some $r \in [1,k-1]$, that a (partial) reverse tableau  $\mathscr R_r$ has been constructed with Red Set being the union of the red multi-sets $J_r^\mathcal C: r\in [1,k-1]$ coming from $\mathscr T^\mathcal C$.

Retain the notation of \ref {3.2.3}. We omit the superscript $\mathcal C$.  Set $C'=C_{r+1}$ and let $h$ be its height. Let $m,m'$ denote the corresponding multiplicities, that is to say $m$ is the cardinality of $J_{r+1}$ as a set and $m'$ is the multiplicity of $j_{r+1}$.

Assume $J_{r+1}$ non-empty, so then $m,m'$ are both $\geq 1$.

By these choices and the rules used in constructing the component tableau $\mathscr T^\mathcal C$ \cite [3.2.2, Rule (1)(b), 3.2.6 Overview]{FJ5} and Remark $1$ of \ref {3.3}, there must be a set $\mathscr P_{r+1}$ of pairs of neighbouring columns of heights running once through $[h-m+1,h+m'-1]$ which are free in the sense that they have not been used in constructing the first $r$ columns of the component tableau $\mathscr T^\mathcal C$. 

Assume that they also have not been used in the construction of $\mathscr R_r$.  (This is established by induction on the construction, without further mention.)

Start with $i=h$.  The assumption that $m\geq 1$ means that $C'$ has a left neighbour $C$ (of course of height $h$) and by assumption none of the pairs in $\mathscr P_{r+1}$, and in particular the pair $C,C'$, have not yet been implemented. Then by \ref {4.4.6} the
 entry $j$ of $R_h\cap C_{r+1}$ which is its lowest entry and is necessarily in black.

 Implementing the pair $C,C'$, recolours $j$ in red and places a black entry $j$ in $R_{h+1}$ in the trapezium $\mathscr T^h_{C,C'}$ as it evolves (cf \ref {4.4.4}) from obtained the rectangle $R^h_{C,C'}$.

This argument may be repeated for $i=h+1, h+2, h+m'-1$, to the black entry which obtains from the previous step.  These are recoloured in red except the very last one.  They define a reverse $j$-string, with $m'$ red with the leftmost one in black.

 Again we may start with $i=h$ and decrease by increments of one to $h-m+1$.  In this the original black
 say $j$ in $C_{r+1}$ is recoloured in red thereby reducing the black height of $C_{r+1}$ successively to $h-1, h-2,\ldots h-m$, whose $m$ entries are eventually  in red.

 This gives the required result.

\end {proof}
\textbf{Remark.}  There is a some limited flexibility in the order of implementation of columns pairs.  Thus we may decrease column heights to $h-m+1$ first and then increase column heights to $h+m'-1$. However we cannot switch column heights ``in the middle''.  This is because we require (in the first case) that the new black entry to lie successively in $R_{i+1}:i=h, h+1,\ldots,h+m'-1$ on the one hand and for the black height of $C_{r+1}$ to decrease successively from $h$ to $h-m+1$.

There two consequences of the above construction which are interesting but are not actually needed.  The reader may reconstruct these results from an earlier unpublished version \cite {FJ6} of this paper, though can  instead remain in blissful ignorance.

\subsubsection {The Appearance of the Amalgamated Column $\hat{C}$ of \cite [4.3.4]{FJ5}}\label {6.1.1}

\

The amalgamated column of \cite [4.3.6]{FJ5} was an ad hoc construction of the first named author and was obtained by empiling Red Sets (in our present language) coming from an ``$i$-string of penetration'' \cite [Definition 4.3.5]{FJ5}.  It is \textit{mandatory} for proving the vanishing of the corresponding invariant on the variety defined by the associated excluded roots.

Mysteriously it does not seem to be needed for the corresponding result (Theorem \ref {5.1}) for the reverse tableau - so what happened?

It turns out that $\mathscr R^{\psi(\mathcal C)}$ incorporated the amalgamated column as in shown in Example $5$ of Section \ref {10}.

In short the reverse tableaux do not need this ``artificial ad hoc'' construction.

 \section {Interchanging Component and Reverse Tableaux}\label {7}

\begin{center} \textbf{}\end{center}

\subsection {A Map from Component to Reverse Tableaux.}\label {7.1}

\

Recall \ref {6.1} the construction of the reverse tableau $\mathscr R^{\psi(\mathcal C)}$.

This gives a map $\Psi: \mathscr T^\mathcal C \mapsto \mathscr R^{\psi(\mathcal C)}$.

\subsection  {The Injectivity of the Map $\Psi$.}\label {7.2}

\begin {lemma}  The map $\Psi$ is injective.
\end {lemma}

\begin {proof}
For each column $C_{r+1};r \in [1,k]$ of $\mathscr D$, it  follows from \ref {6.1} that the reverse tableau determines both the red subset of the $C_{r+1}$ of the component tableau from which it was derived and the multiplicity of its lowest entry $j$ as the number of appearances of a red $j$ in the reverse tableau.  Then by the injectivity of the red map (\ref {3.3}), it follows that $\Psi$ is injective.

\end {proof}

\subsection  {The Bijectivity of the Map $\Psi$.}\label {7.3}

Bijectivity of $\Psi$ fails because there can be many different reverse tableaux with the same Red Set (see the right hand column of Fig. $14$, for example).  However by the Remark to Theorem \ref {9.3}, reverse tableaux with the same Red Set define the same component of $\mathscr N$ via their excluded roots.

Identifying such reverse tableaux makes $\Psi$ bujective.

\section {Some Notation, Notions and Further Comments }\label {8}

\subsection {Notation}\label {8.1}

For a given tableau (to be written in parentheses), let $X$ be the set of excluded root vectors, $Y$ the set of root vectors labelled by a $\ast$ and $S$ the set of root vectors labelled by a $1$, it defines.

We require that

\

$(i)$. $X \supset Y$ and we set $Z=X\setminus Y$.

\

$(ii)$. $S\cap X=\phi$.

\

$(iii)$. As vectors in \textbf{M}, the elements of $S$ lie in distinct rows and distinct columns\footnote{Pictured as rooks on a chessboard  none cover any other.} .

\

As noted in  \cite [6.3(i)]{FJ5}, property $(iii)$  makes calculations of tangent space dimension based on them much easier.  In particular it would have been no advantage if $S$ were (to be made) bigger, if $(iii)$ is lost in the process.

\

When $X,Y,Z,S$ refer to a particular tableau, be it component or reverse, the latter is placed in parentheses.

For a component tableau, properties $(i,ii,iii)$ above are established in \cite [Lemma 4.1.5]{FJ5}, \cite [Prop. 4.2.5]{FJ5} and \cite [Cor. 3.2.6(ii)]{FJ5}, respectively.  For a reverse tableau, $(i)$ is established in Lemma
 \ref {8.2}, $(ii)$ in Lemma \ref {8.3}.  Moreover $(iii)$ for the reverse tableau obtains from $(iii)$ for the component tableau since by imposing

$$Y:=Y(\mathscr R^{\psi(\mathcal C)})=Y(\mathscr T^\mathcal C), \quad S:=S(\mathscr R^{\psi(\mathcal C)})=S(\mathscr T^\mathcal C),\eqno (*)$$
the consequences of which we examine in \ref {8.2}, \ref {8.3} below.


\subsubsection {Reverse strings}\label {8.1.1}

Recall the notion of a reverse $i$-string \ref {4.3}.%

 The left (right) end point of a reverse $j$-string is black (resp. red). On the other hand if  $m'=0$ then the reverse $j$-string degenerates to a single entry with just one black label.

\subsection {Lines with label $\ast$ in the Reverse Tableaux}\label {8.2}

\

Let $\ell_{i,j}$ be a line with label $\ast$ in the component tableau $\mathscr T^\mathcal C$.

The imposed identity in Eq. $(*)$ of Section \ref {8.1}, means that we can choose a line $\ell_{i,j}$ with label $\ast$ in the corresponding reverse tableau $\mathscr R^{\psi(\mathcal C)}$ with the same numerical end points, namely $i,j$.

Of course this would be rather optimistic had not the Red Sets coincided in the reverse and component tableaux. Indeed in the latter,  $j$ is the unique entry in red (with $m'$ its multiplicity, \textit{necessarily} $\geq 1$) and $i$ the leftmost entry (which is in black). (This condition implies that both  lie in $\mathscr T^\mathcal C$ rather than just in $\mathscr T^\mathcal C(\infty)$.)


Thus we take $j\in \mathscr R^{\psi(\mathcal C)}$ to be the unique black entry with numerical value $j$.

On the other hand the numerical entries of $i$ in the corresponding reverse tableau $\mathscr R^{\psi(\mathcal C)}$ form a unique reverse $i$-string (\ref {4.3}) and we choose $i$ to be its unique rightmost entry. Recalling that Red Sets and so multiplicities coincide, so then this reverse string has $m'+1$ elements.  Thus its rightmost entry is in red!

Recall the notation introduced in Remark $1$ of \ref {4.3.6}.

\begin {lemma}  With these choices $\ell_{i,j}^{\psi(\mathcal C)} \in X(\mathscr R^{\psi(\mathcal C)})$.
\end {lemma}

\begin {proof}


For each $i \in \mathscr T^\mathcal C(\infty)$ joined by a vertical line (labelled by a $\ast$) to some $j \in J^\mathcal C_{r+1}$, one has $i<j$, because $i$ had been shifted strictly from its original place in $\mathscr T^\mathcal C$ in turn strictly to the left of $C_{r+1}$.

On the other hand in the construction of the reverse tableau \textit {because $m'\geq 1$}, $j$ is placed under $i$ in the same column of $\mathscr T^\mathcal C$, which then may be skewed to the left by subsequent operations.


Thus the line starting at $j\in B(\mathscr R)$ to $i$ goes rightwards (not necessarily strictly).  Since $i<j$ this forces $x_{i,j}$ to be an excluded root for the reverse tableau.

%
%
 \end {proof}


%

%

\subsection {Lines with label $1$ in the Reverse Tableaux}\label {8.3}

\

Let $\ell_{i,j}$ be a line with label $1$ in the component tableau.

For the corresponding line $\ell_{i,j}$ with label $1$ in the reverse tableau $\mathscr R^{\psi(\mathcal C)}$, Eq. $(*)$ of Section \ref {8.1} requires that their numerical labels be preserved.

 As in \ref {8.2}, we similarly resolve the ambiguity of choice in the reverse tableau.



Thus in the reverse tableaux, we take $j$ to be the unique entry in black and $i$ to be the rightmost entry which may be either in black (for $m'=0$) and in red otherwise.

In a reverse tableau the unique black element with label $j$ lies $m'$ rows below the unique rightmost red element in the reverse $j$ string. (This comes from combining ``Substitution'' which introduces a new entry one row below and ``Shifting of Partial of columns'' which leaves rows unchanged \ref {4.3}.

Yet $\ell_{i,j}$ is upward going  by $m'$ rows in $\mathscr T^\mathcal C(\infty)$ and this exactly cancels the black element $j$ being $m'$ rows below its rightmost entry.

In particular if $i\in \mathscr T$ appears in the same row as it does in $C_r(\infty)$, then the line $\ell_{i,j}$ in the reverse tableau is horizontal.

   More generally $i$ lies in some column to the left of $C_r$, some $u\geq 0$ rows higher up than it appears in $C_r(\infty)$.

   \

   Observe, by say \cite [Fig. 1 of Sect. 8]{FJ5}, that $u$ is also the number of vertical lines with label $\ast$ from $i$ to appropriate columns in $\mathscr T^\mathcal C(\infty)$  to boxes in $\mathscr T^\mathcal C$ with entries necessarily strictly greater than $i$.

   \

   We conclude that in terms of \textbf{M} one has

 \

 $(*)$.  Suppose in row $i$ of $\textbf{M}$, that there is an entry labelled by $1$  has $u$ entries on its left labelled by a $\ast$.  Then in the reverse tableau the right going line $\ell_{i,j}:i<j$ with label $1$, described above goes down by $u\geq 0$ rows.

 \

$(**)$. With the exception described in $(*)$ the left going lines labelled by a $1$ in a reverse tableau are horizontal.  Without exception they are all down going.


\begin {lemma} $S(\mathscr R^{\psi(\mathcal C)})\cap X(\mathscr R^{\psi(\mathcal C)})=\phi.$
\end {lemma}

\begin {proof} Indeed by $(**)$, the lines with label $1$ are right and down-going so by \ref {8.1} are not excluded in $\mathscr R^{\psi(\mathcal C)}$.
\end {proof}
\textbf{Example 8.}  Consider the composition $(2,1,1,2,2)$ and the composition with red multi-set $4,6,8$.   In row $3$ of \textbf{M} the entry $1$ coming from the line $\ell_{3,8}$ follows two entries of $\ast$ (in columns $4,6$) and so in the reverse tableau going rightwards descends by two rows.   See blue line in Sect. \ref {10}, Figure $8$.

\

 \textbf{Example 9}.  Consider the composition  $(3,2,1,3,2,1,2)$ and the component with red multi-set $8,9,12,12$.  The line $\ell_{9,12}$ goes up by $m=2$ rows in the component tableau whilst the new black $12$ goes down by $m=2$ rows in the reverse tableau.  Thus the line $\ell_{9,12}$ labelled by a $1$ is horizontal in the reverse tableau. See blue line in Sect. \ref {10}, Figure $9$.

 \section {Covering}\label {6}

\subsection {Components}\label {9.1}

\subsubsection {Coincidences}\label {9.1.1}

 let $\mathfrak u^\mathcal C$ denote the subspace of $\mathfrak m$ spanned by the complement to the excluded roots in $X(\mathscr T^\mathcal C)$.

 Set  $\mathscr C:=\overline {B.\mathfrak u^\mathcal C}$, which is a closed irreducible subvariety of $\mathfrak m$. It is shown in \cite [Prop. 6.3]{FJ5} through ``covering'' that
 $$ \dim \mathscr C \geq \dim \mathfrak m - \textbf{g}. \eqno {(E)}$$

 Then in \cite [Thm. $1$]{FJ5} that  $\mathscr C$ is a component of $\mathscr N$ and equality holds in $(E)$.


  Similarly let $\mathfrak u^{\psi(\mathcal C)}$ denote the subspace of $\mathfrak m$ spanned by the complement to the excluded roots in $X(\mathscr R^{\psi(\mathcal C)})$.

  We wish to show that $\mathscr C':=\overline {B.\mathfrak u^{\psi(\mathcal C)}}$, which is a closed irreducible subvariety of $\mathfrak m$, is a component of $\mathscr N$ in $\mathfrak m$ of codimension \textbf{g} and furthermore that it coincides with $\mathscr C$.

  Notice that this last assertion will also imply that two reverse tableaux with the same Red Set define the same component via excluded roots.

  \subsubsection {Use of a Weierstrass Section.}\label {9.1.2}

  The first step towards the above goal is exactly as in \cite [6.1]{FJ5}.

  Thus we assume that the inequality
  $$ \dim \mathscr C' \geq \dim \mathfrak m - \textbf{g}, \eqno {(E')}$$
  holds.

  If $\mathscr C' $ is not a component of $\mathscr N$, then it is strictly contained in some component $\hat{\mathscr C}$.  By $(E)$ this must have dimension strictly greater than $\dim \mathfrak m - \textbf{g}$.  Moreover it is conical and its projectivisation has dimension $\geq \dim \mathfrak m - \textbf{g}$.

  Recalling \cite [Thm. 5.2.6]{FJ5} let $e_\mathcal C+V_\mathscr C$ be the Weierstrass section associated to $\mathscr T^\mathcal C$.

  Let $\mathscr V^e$ be the set of all non-zero multiples of $e_\mathcal C+v: v\in V_\mathscr C\setminus \{0\}$, which is again conical.  Its projectisation has dimension \textbf{g}-1.

   Then, as noted in \cite [Thm. 6.1]{FJ5}, it follows from the intersection theory in projective space that $\hat{\mathscr C}$ must contain an element of the form $e_\mathcal C+v: v\in V_\mathscr C\setminus \{0\}$ and hence that there is a Benlolo-Sanderson invariant that does not vanish on $\hat{\mathscr C}$. This contradicts Theorem \ref {5.1}, thereby proving that $\mathscr C'$ is indeed a component of $\mathscr N$ of dimension $\dim \mathfrak m-\textbf{g}$, all the time assuming that $(E')$ holds.

   \subsubsection {Coincidence of $\mathscr C, \mathscr C'$.}\label {9.1.3}

   We may also give a dimension argument to obtain the coincidence of $\mathscr C, \mathscr C'$.  Indeed let $\mathfrak u^{\mathcal C \cup \mathcal \psi(\mathcal C)}$ denote the subspace of $\mathfrak m$ of roots not excluded either in  $X(\mathscr T^\mathcal C)$ nor in $X(\mathscr R^{\psi(\mathcal C)})$ and set $\mathscr C''=\overline {B.\mathfrak u^{\mathcal C \cup \mathcal \psi(\mathcal C)}}$.  This contains $\mathscr C, \mathscr C'$.  Thus if we can show that
     $$ \dim \mathscr C'' \geq \dim \mathfrak m - \textbf{g}, \eqno {(E'')}$$
     then all three varieties must coincide.  In particular $\mathscr C = \mathscr C'$, as required.

     As in the proof $(E)$, the proof of $(E'),(E'')$ is by covering and then using \cite [Prop. 6.3]{FJ5}.

  \subsection {The Combinatorics of Covering}\label {9.2}

  \

   Let $\mathfrak n$ be the Lie algebra of the nilradical of $\mathfrak b$.

  \subsubsection {Definition of Covering}\label {9.2.1}

  \

  We recall the notion of covering \cite [6.2.2]{FJ5}.  Fix a row of \textbf{M} and let $x,y$ be distinct co-ordinate vectors in that row.  We say that $x$ covers $y$ if $x$ lies (strictly) to the left of $y$.  This notion arises since then $y \in \mathfrak b.x$, that is to say that $y$ lies in the tangent space to $B$ at $x$.

The application of the above simple fact is more complicated since we shall not use a single co-ordinate vector, but rather a sum $e_\mathcal C =\sum_{x\in S} x$,  recalling the notation of \ref {8.1}, Eq. $(*)$ that $S=S(\mathscr R^{\psi(\mathcal C)})=S(\mathscr T^\mathcal C)$.

Nevertheless in view of \ref {8.1}(iii), an induction argument (given in the proof of \cite [Prop. 6.3(i)]{FJ5}) that if $S$ covers $Z=X\setminus Y$  then $\mathfrak u + \mathfrak n.e +Y=\mathfrak m$, where one's tableau of choice is written in parenthesis and $\mathfrak u$ is the subspace of $\mathfrak m$ spanned by the roots in the corresponding tableau which are not excluded.

\subsubsection {The Key Proposition}\label {9.2.2}

Here we shall write simply $\mathscr T^\mathcal C$ (resp. $\mathscr R^{\psi(\mathcal C)}$) as $\mathscr T$ (resp. $\mathscr R$).

\begin {prop} $x_{i,j}$ is covered in $\mathscr R$ if and only if $i$ is stopped in $\mathscr T$ before it reaches $C$.
\end {prop}

The proof will be given in the subsections below.

\subsubsection {The Stopping Entry}\label {1.3.1}  If $i$ is stopped in $\mathscr T(\infty)$ before it reaches $C$, then there is an entry $k$ to the left (not necessarily strictly) of $C$ such that $x_{i,k}$ is labelled by a $1$ in \textbf{M}.  (Loosely speaking $k$ is the label that stops $i$, that is the ``stopping entry''.) Then to show that the excluded root is covered (by $\ell_{i,k}$) it suffices to show that $k<j$.

The inequality $k<j$ will be shown in \ref {9.2.4}.  Let us examine its consequence.

If $i$ is not stopped before it reaches $C$, then $i$ enters $C(\infty)\setminus C$ in $\mathscr T(\infty)$, and there is a vertical line $\ell_{i,j}$ in $\mathscr T(\infty)$ labelled by a $\ast$.  Such a line is not required to be covered and since their total number is just \textbf{g} by \cite [Lemma 3.2.5]{FJ5}, it follows
that it will not be covered.

\textit{In other words in the proposition ``if'' also implies ``only if''.} via the last part of \ref {1.2}, that is to say by comparing the total number of lines with label $\ast$ to the number of Benlolo-Sanderson invariant generators.

\subsubsection {Proof of the Inequality $k<j$}\label {9.2.4}

\

If $i$ is stopped strictly before it reaches $C$, then the numbering \ref {2.5}) in $\mathscr T$ implies automatically that $k<j$.

In the case that $i$ is stopped exactly at $C$, then $j$ is chosen inductively (\cite [3.2.6]{FJ5}) to be the highest available element in right adjacent column $C^+$ to $C$ as its highest available element, whose entry is not already the end-point of a line in $C(\infty)$.  The resulting pattern of lines labelled by $1$ or by $\ast$ is described in \cite [{3.2.6}, Overview]{FJ5} and illustrated in \cite [Sect. $9$, Fig $1$]{FJ5}.

 From the latter it is particularly clear that $k<j$.  Indeed in its right hand column $k$ must be one of the end-points of one of the diagonal black lines (or the single right of which there is only one - see caption to figure!) and up-going dashed/dotted (blue) line, whilst $j$ is the end-point of the right and down going dashed (green) lines.  (The stated colour scheme being valid depends on the combined obligeance of the journal editor, not being colour-blind and the reader's printer. Otherwise lines can be distinguished by being in bold, having dashes etc.

 \subsubsection {Conclusion} {\label {9.2.5}

 \

 \begin {cor}

 \

 $(i)$.  The set $S(\mathscr R^{\psi(\mathcal C)})$ does not cover  $Y(\mathscr R^{\psi(\mathcal C)})$.

  \

  $(ii)$. The set $S(\mathscr R^{\psi(\mathcal C)})$ covers the roots of $Z(\mathscr R^{\psi(\mathcal C)}):= X(\mathscr R^{\psi(\mathcal C)})\setminus Y(\mathscr R^{\psi(\mathcal C)})$.
 \end {cor}

\subsection {Associated varieties}\label {9.3}

We recall that a similar result to \ref {9.2.5}  was established for the component tableau $\mathscr T^\mathcal C$ in \cite [6.2.4-6.2.7] {FJ5}, except that the proof was \textit{considerably more painful}.  Through \cite [Prop. 6.3]{FJ5} we concluded from \cite [Thm. 6.1]{FJ5}  that $\mathscr C$ is a component of $\mathscr N$ of codimension \textbf{g} in $\mathfrak m$. Similarly (as explained) above we obtain from Corollary \ref {9.2.5} that $\mathscr C'$ is a component of $\mathscr N$ of codimension \textbf{g} in $\mathfrak m$.

Finally since the common set $S$ covers both the roots \textit{not} excluded from $X(\mathscr T)$ and \textit{not} excluded from $X(\mathscr R)$ - just as a trivial consequence of the definition of being excluded.  This gives $(E')$ and may be summarized in the
\begin {thm}  $\overline {B.\mathfrak u^\mathcal C} = \overline {B.\mathfrak u^{\psi(\mathcal C)}}$.
\end {thm}
\textbf{Remark.}  Thus the variety associated to a reverse tableau via its excluded roots only depends on its Red Set.  By the injectivity of the component map \cite [Sect. 7]{FJ5} and the injectivity of the red map \ref {3.3}, these are pairwise distinct.  In particular if we identify reverse tableaux which admit the same red set, then the map $\Psi$ of \ref {7.2} is bijective.

\subsection {A Shortcut to Covering in the Component Tableau}\label {9.4}

The proof of covering and of  \ref {8.1}(ii) for the reverse tableau $\mathscr R$ is much easier than that for the corresponding component tableau $\mathscr T$ - see \ref {8.1}.  Thus it is relevant to point out that the former implies the latter when $X(\mathscr R)\supset X(\mathscr T)$.  In the preliminary version of the present work \cite {FJ6},  ``redressing'' is used to show that for any red set there is the reverse tableau satisfying this inclusion using extreme shifting \cite [Lemma 4.5 and Prop. 6.3.1]{FJ6}. No doubt this can also be proved using ``enabling''.

\subsection {A Glimpse into Factorisation}\label {9.5}

To whet the reader's appetite we illustrate the factorisation procedure outlined in \ref {3.6}, through a flow chart in Example 14, describing the two choices of complete sequences to recover the expected reverse tableaux.    The general case will be studied in ``Enabling''. This example shows that the enabling proposition is really needed - \textit{redressing is not enough}. Again in this example one checks that the enabling proposition at each stage holds by inspection.

As outlined in \ref {3.6}, we must consider enough sequences of all invariant generators so that successive linearisation of the first $i$ invariants and factorisation of the next gives all the components we already know exist. Every such sequence corresponds to a complete sequence of pairs of neighbouring columns. Here by Krull's theorem, the factorisation of the $i^{th}$ invariant given by the reverse tableau on implementing the next pair of neighbouring columns $(C_i, C_{i'})$, through the enabling proposition.  These reverse tableaux are defined by their red multi-set indicated as a subscript.  It is the same red multi-set as for some component tableaux \cite {FJ7}.

\

Consider now the composition $(1,2,3,1,1,3,2)$.  The component tableau gives $6$ red multi-sets (Example $1$) namely
$(7,7,7,8), (7,8,8,8), (7,7,8,11), (7,8,8,11),(7,8,11,10), (7,8,11,13)$.
All contain the Red Set  $(7,8)$, so we start at the corresponding (partial) reverse tableau, which we designate by a subscript, that is to say as $\mathscr R_{7,8}$ in the present case.


Implement the pair $(C_2,C_7)$ (of neighbouring columns of $\mathscr D$) on $\mathscr R_{7,8}$.  Since  $\mathscr R_{7,8}$ has $4$ columns of height $2$ and since the enabling proposition holds by inspection we obtain three reverse tableaux
$\mathscr R_{7,7,8}, \mathscr R_{7,8,8}, \mathscr R_{7,8,13}$.

Now implement the pair $(C_3,C_6)$ to each such reverse tableau.  Each has three columns
of height $3$. The first gives the two reverse tableaux $\mathscr R_{7,7,7,8},\mathscr R_{7,7,8,1}$.  The second gives $\mathscr R_{7,8,8,8},\mathscr R_{7,8,8,11}$. The third just gives \textit{only} $\mathscr R_{7,8,11,13}$,
 since the
\textit{the original black} entry lies in $C_7$ and so is not
surrounded by $(C_3,C_6)$ and is excluded by the rule in \ref {4.1.1}.

This rule is equivalent to the rule that implementing a neighbouring pair $C_i,C_{i'}$ of height $s$ one cannot drop an entry strictly to the right of a $C_{i'}$ from row $R_s$ to row $R_{s+1}$.
It would have given a reverse tableau $\mathscr R_{7,8,13,13}$.  Again there is no composition tableau with this as its red multi-set.

  Again it is also obvious that the invariant $I^3_{C_3,C_6}$ cannot factor since the  only excluded root is $x_{7,8}$ in the sub-matrix of \textbf{M} in $\mathfrak m$ bounded by the Levi blocks of size three.  This fact was the origin by the rule in \ref {4.4.4}.

Now implement the pair $(C_3,C_6)$. Since $\mathscr R_{7,8}$ has two columns of height $3$,
we just get one reverse tableau $\mathscr R_{7,8,11}$. This has five columns of \textit{black} height $2$. Since
$C_2^-=C_1$  we obtain
four reverse tableaux, namely $\mathscr R_{7,7,8,11},\mathscr R_{7,8,8,11}.\mathscr R_{7,8,10,11},\mathscr R_{7,8,11,13}$.

One reverse tableau in the first sequence is missing, two in the second sequence. Again the reverse tableau $\mathscr R_{7,8,8,11}$ differs slightly from that previously
but by the Remark to Theorem \ref {9.3}, they nevertheless define the same variety.

Finally to obtain the surjectivity of the component map we must eliminate the elephant in the room. Fortuitously, as noted in \ref {3.6}, this reduces to showing that linearisation does not eliminate any components.


Perhaps the reader's head is spinning from all this intermingling of tableaux; but pity rather the second author who gamely tried to make the daring ideas of the first author fit for human consumption.

\subsection {Index of Notation and of Notions}\label {6.8}

Notations and notions frequently used are given below in the Section (excepting the Abstract) where they are first fully defined.

\

\textbf{Notations}

\ref {1.2}.  \ \ \ \ $\mathbb C,[m,n],]m,n],[m,n[,]m,n[$.

\ref {2.1}.  \ \ \ \  $\mathfrak {gl}(n),\textbf{M}_n,\mathfrak {sl}(n), \mathfrak b, \mathfrak p, \mathfrak r,\textbf{B}_i,\mathfrak m$.

\ref {2.2}. \ \ \ \  $\mathscr D, C_i,\textbf{C}_i,c_i,x_{i,j},x^\ast_{i,j}, R_i,R^u,R^{>u},b_{i,j}$.

\ref {2.3}. \ \ \ \ $[C,C'],]C,C'],[C,C'[,]C,C'[,\textbf{g}, C^{>u},I^s_{C,C'},\mathscr I$.

\ref {2.4}. \ \ \ \ $R^s_{C,C'}, \ ^lR^s_{C,C'}$.

\ref {2.5}. \ \ \ \ $\mathscr T, \ell_{i,j}^\mathcal C$.

\ref {2.7}. \ \ \ \ $I^s_{C,C'},\mathcal I, \mathscr N$.

\ref {3.1}. \ \ \ \ $\mathscr T^\mathcal C, S(\mathscr T^\mathcal C), Y(\mathscr T^\mathcal C)$.

\ref {3.2}. \ \ \ \ $X(\mathscr T^\mathcal C), \mathfrak u^\mathcal C, \mathscr C$.

\ref {3.2.1}. \ \  $\mathscr T^\mathcal C(\infty)$.

\ref {3.2.2}. \ \  $\mathscr B^s_{C,C'}$.

\ref {3.2.3}. \ \  $S(i), J^\mathcal C_{r+1}$.

\ref {3.2.4}. \ \  $m_{r+1}^\mathcal C, j_{r+1}, m'_{r+1}, h_{r+1} R(\mathscr T^\mathcal C)$.

\ref {3.2.5}. \ \ $\ell_{i,j}^\mathcal C$.


\ref {3.5}. \ \ \ \ $\textit{p},\mathscr P,\mathscr P_t$.


\ref {4.1.1}. \ \ $\mathscr I_i, C^-$.

\ref {4.2}. \ \ \ \ $S_i$.

\ref {4.3.3}. \ \   $R(\mathscr R)$.

\ref {4.3.4}. \ \  $X(\mathscr R^\mathcal E),  X(\mathscr T^{\mathcal C}_{i,j})$.

\ref {4.4.1}. \ \ $B^i$.

\ref {4.4.2}. \ \ $B'^i$.

\ref {6.1}. \ \ \ \ \  $\mathscr R^{(\psi(\mathcal C)}$.

\ref {7.1}. \ \ \ \ \ $\Psi$.

\ref {8.1}. \ \ \ \ \  $X(\cdot),Y(\cdot), Z(\cdot),S(\cdot)$.

\ref {8.2}. \ \ \ \ \ $\ell^{\psi(\mathcal C)}$.

\ref {9.1.1}. \ \ \ $\mathfrak u^\mathcal C, \mathfrak u^{\psi(\mathcal C)}$.

\

\textbf{Notions}.

\

\ref {2.1}. \ \ \ \ Blocks.

\ref {2.2}. \ \ \ \ Column blocks, Rows, boxes.

\ref {2.3}. \ \ \ \ Adjacent and Neighbouring Columns, Intervals, Lower parts of columns.

\ref {2.4}. \ \ \ \ Rectangles.

\ref {2.6}. \ \ \ \ Surrounding Columns.

\ref {2.7}. \ \ \ \ Benlolo-Sanderson Invariant (BS).

\ref {3.1}. \ \ \ \ Line decoration, Component Tableaux.

\ref {3.2}. \ \ \ \ Excluded Roots, The component map.

\ref {3.2.1}. \ \ Semi-standard tableau.

\ref {3.2.2}. \ \ Batches.

\ref {3.2.3}. \ \ Strings, neutral lines, free pairs, stopped strings.

\ref {3.2.4}. \ \ The Red Set.

\ref {3.2.5}. \ \ Red lines.

\ref {3.3}. \ \ \ \ The red map, Canonical component.

\ref {3.5}. \ \ \ \ Complete sequences of neighbouring columns.

\ref {3.6}. \ \ \ \ The Elephant.

\ref {4.1}.\ \ \ \ \ \   Reverse tableau, black height.

\ref {4.1.1}. \ \  The left boundary column.

\ref 4.2}. \ \ \ \ \ The Enabling Proposition.

\ref {4.3.1}. \ \ \ Substitution, Shifting of Partial columns, Skewing to the left.

\ref {4.3.2}. \ \ \  Extreme shifting.

\ref {4.3.3}. \ \ \  Reverse $j$-strings, Red Set of the reverse tableau.

\ref {4.4.1}. \ \ \  The left boundary.

\ref {4.4.2}. \ \ \  The right boundary.

\ref {4.4.3}. \ \ \  The trapezium.

%
%

\ref {4.3.4}. \ \ \  Standard with multiplicities.

\ref {4.3.5}. \ \ \  Excluded roots.

\ref {6.1.1}. \ \ \  Amalgamated column.

\ref {9.2.1}. \ \ \ Covering.

\ref {9.5}. \ \ \ \ \ \ Factorization.

\section {Illustrated Examples.} \label {10}



\textbf{Example 1}.

Consider the composition $(1,2,3,1,1,3,2)$.  The six possible component tableaux $\mathscr T^\mathcal C(\infty)$ are computed below using the rules given in \cite [3.2.2]{FJ5}.  Then Red Sets $J^\mathcal C_{r+1}$, as defined in \ref {3.2.3}, are illustrated in red. In addition the multiplicity of its lowest element $j_{r+1}$, defined in \ref {3.2.4}, can be read off.  This leads to $6$ Red Sets
\[
(7,7,7,8),\ (7,8,8,8),\ (7,7,8,11),\ (7,8,8,11),\ (7,8,11,10),\ (7,8,11,13).
\]
one for each of the component tableaux. The first two pairs coincide as sets.
\begin{figure}[H]
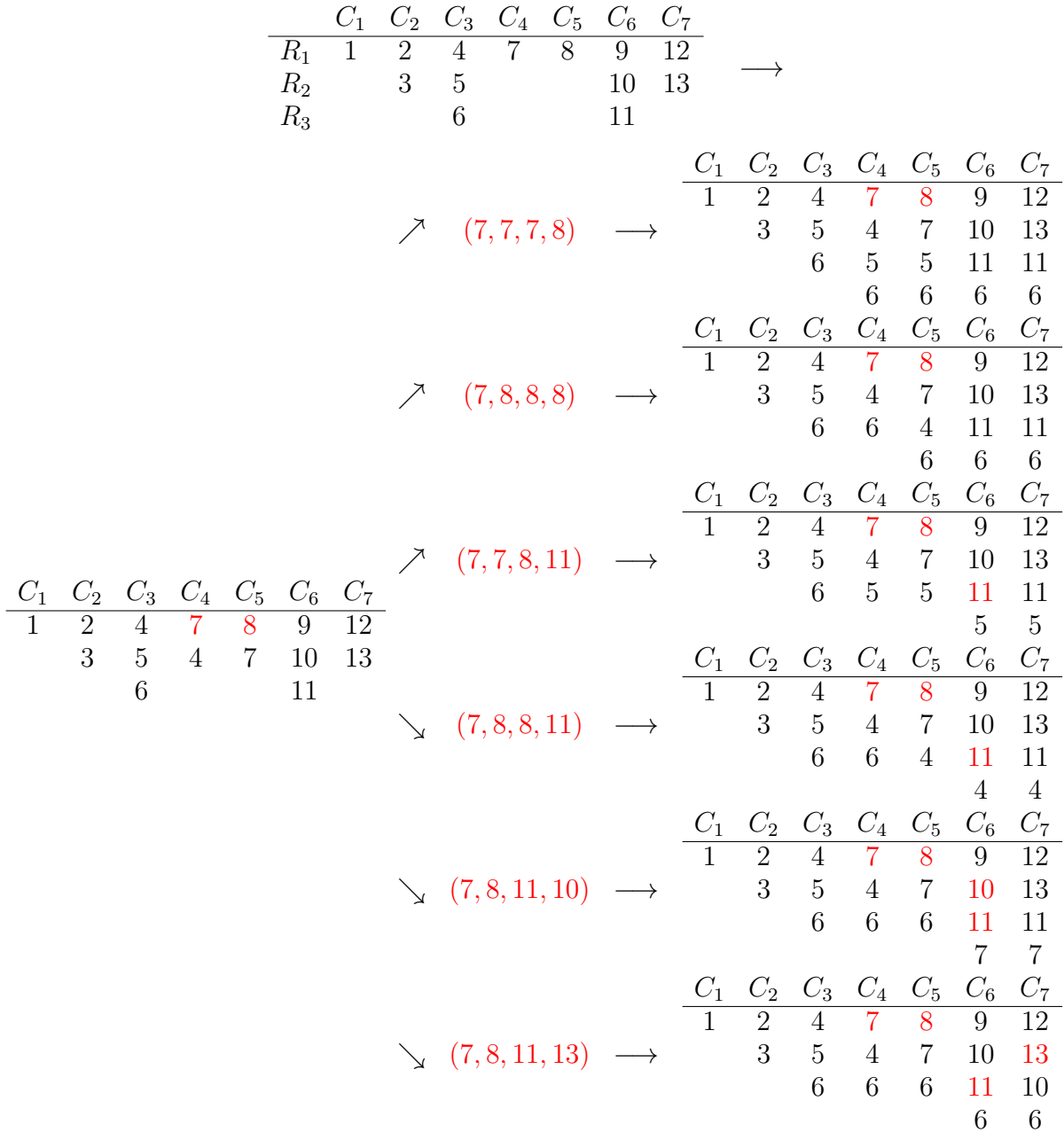

 \[\begin{array}{ ccccccccc}
 &C_1 & C_2 & C_3 & C_4 & C_5 & C_6 & C_7 \\ \hline
  R_1 &1&2&4&7&8&9&12\\
    R_2& &3&5& &&10&13\\
    R_3&& & 6& &&11&
  \end{array} \quad \longrightarrow \quad
  \]
  \[\begin{array}{ cccccccc}
 C_1 & C_2 & C_3 & C_4 & C_5 & C_6 & C_7 \\ \hline
1&2&4& \textcolor{red}{7}& \textcolor{red}{8}&9&12\\
 &3&5&4 &7&10&13\\
& & 6& &&11&
  \end{array} \begin{array}{cccc}
  \nearrow &\textcolor{red}{(7,7,7,8)}&\longrightarrow &\begin{array}{ cccccccc}
 C_1 & C_2 & C_3 & C_4 & C_5 & C_6 & C_7 \\ \hline
1&2&4& \textcolor{red}{7}& \textcolor{red}{8}&9&12\\
 &3&5&4 &7&10&13\\
& & 6& 5&5&11&11\\
& & & 6&6&6&6
  \end{array}\\
  \nearrow& \textcolor{red}{(7,8,8,8)}&\longrightarrow &\begin{array}{ cccccccc}
 C_1 & C_2 & C_3 & C_4 & C_5 & C_6 & C_7 \\ \hline
1&2&4& \textcolor{red}{7}& \textcolor{red}{8}&9&12\\
 &3&5&4 &7&10&13\\
& & 6& 6&4&11&11\\
& & & &6&6&6
  \end{array}\\
  \nearrow&\textcolor{red}{(7,7,8,11)}&\longrightarrow &\begin{array}{ cccccccc}
 C_1 & C_2 & C_3 & C_4 & C_5 & C_6 & C_7 \\ \hline
1&2&4& \textcolor{red}{7}& \textcolor{red}{8}&9&12\\
 &3&5&4 &7&10&13\\
& & 6& 5&5& \textcolor{red}{11}&11\\
& & & && 5&5\\
  \end{array}\\
    \searrow&\textcolor{red}{(7,8,8,11) }&\longrightarrow &\begin{array}{ cccccccc}
 C_1 & C_2 & C_3 & C_4 & C_5 & C_6 & C_7 \\ \hline
1&2&4& \textcolor{red}{7}& \textcolor{red}{8}&9&12\\
 &3&5&4 &7&10&13\\
& & 6& 6&4&\textcolor{red}{11}&11\\
& & & &&4&4
  \end{array}\\
\searrow&\textcolor{red}{ (7,8,11,10)}&\longrightarrow &\begin{array}{ cccccccc}
 C_1 & C_2 & C_3 & C_4 & C_5 & C_6 & C_7 \\ \hline
1&2&4& \textcolor{red}{7}& \textcolor{red}{8}&9&12\\
 &3&5&4 &7&\textcolor{red}{10}&13\\
& & 6& 6&6&\textcolor{red}{11}&11\\
& & & &&7&7\\
  \end{array}\\
\searrow&\textcolor{red}{(7,8,11,13) }&\longrightarrow &\begin{array}{ cccccccc}
 C_1 & C_2 & C_3 & C_4 & C_5 & C_6 & C_7 \\ \hline
1&2&4& \textcolor{red}{7}& \textcolor{red}{8}&9&12\\
 &3&5&4 &7&10&\textcolor{red}{13}\\
& & 6& 6&6&\textcolor{red}{11}&10\\
& & & &&6&6\\
  \end{array} \\
  \end{array}\]
  \caption{ The component tableaux lead to $6$ distinct and explicitly describable components of the nilfibre (in this case) which we wish eventually to prove is a complete set. In the last example every column has just one red entry and with multiplicity one.  This property characterizes the ``canonical component'',  \cite [1.4] {FJ4} being implicit in the construction of the composition tableau \cite [4.7]{FJ4}.  In the canonical component every string descends by at most row on passing through a given adjacent pair of columns and no two different strings enter the same column.}
  \end{figure}

\textbf{Example $2$}.
Consider the composition $(1,2,1,2)$

\begin{figure}[H]
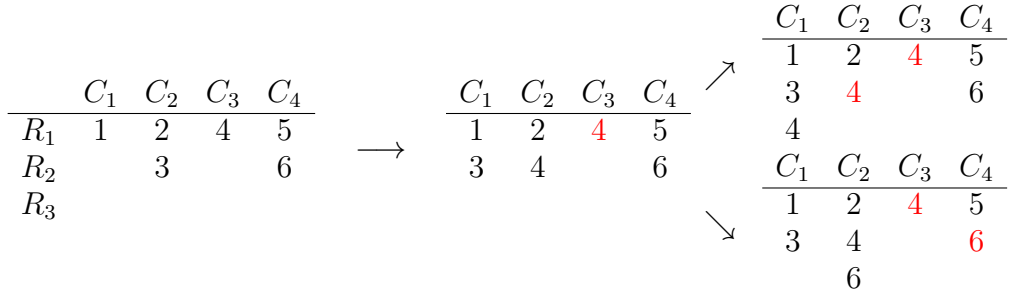

\centering
\[\begin{array}{ccccc}
 & C_1 & C_2 & C_3 & C_4 \\ \hline
R_1 & 1 & 2 & 4 &5  \\
R_2 &   & 3 &  &  6 \\
R_3 &   &   &  &   \\
\end{array}\quad \longrightarrow \quad  \begin{array}{cccc}
 C_1 & C_2 & C_3 & C_4 \\ \hline
 1 & 2 & \textcolor{red}{4} &5  \\
   3& 4 &  &  6 \\
   &   &  &   \\
\end{array} \begin{array}{cc}\nearrow& \begin{array}{cccc}
 C_1 & C_2 & C_3 & C_4 \\ \hline
 1 & 2 & \textcolor{red}{4} &5  \\
   3&  \textcolor{red}{4}  &  &  6 \\
   4 &   &  &   \\
\end{array}\\
\searrow& \begin{array}{cccc}
 C_1 & C_2 & C_3 & C_4 \\ \hline
 1 & 2 & \textcolor{red}{4} &5  \\
   3& 4 &  & \textcolor{red} 6 \\
   &   6&  &   \\
\end{array}

\end{array}\]
\caption{ Implement $C_1,C_3$.  This puts a black $4$ in $C_2$ pushing the black $3$ into $R_2\cap C_1$.  Then for the pair $C_2,C_4$ being implemented, one has $C_2^-=C_1$ of height $2$ with lowest entry in black and hence of height and black height $2$, 
This example gives eventually two red multi-sets $(4,4);(4,6)$.}
\end{figure}
 Example $2_1$. Components tableaux for the composition $(1,2,1,2)$ above: 
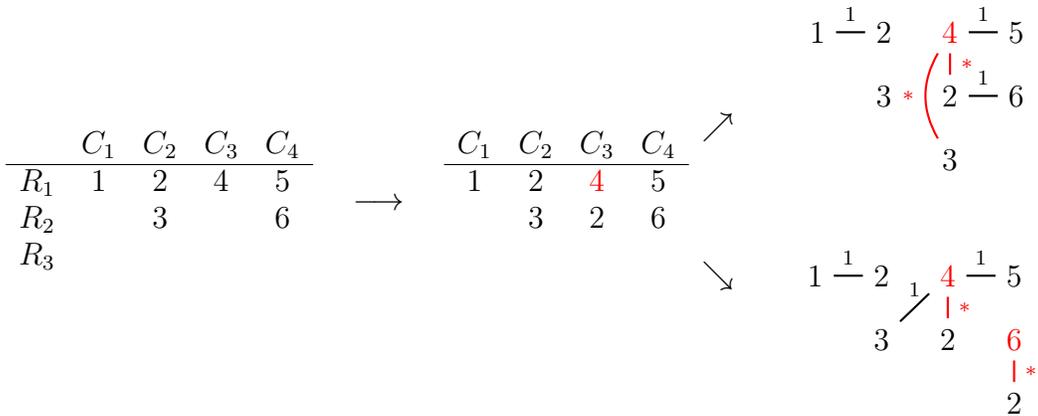
\begin{figure}[H]
\centering
\[\begin{array}{ccccc}
 & C_1 & C_2 & C_3 & C_4 \\ \hline
R_1 & 1 & 2 & 4 &5  \\
R_2 &   & 3 &  &  6 \\
R_3 &   &   &  &   \\
\end{array}\quad \longrightarrow \quad  \begin{array}{cccc}
 C_1 & C_2 & C_3 & C_4 \\ \hline
 1 & 2 & \textcolor{red}{4} &5  \\
   &  3& 2 &  6 \\
   &   &  &   \\
\end{array}\begin{array}{c} \nearrow\\ 
\\
\\
\\
\searrow\end{array}
\begin{array}{cc}
& \begin{tikzpicture}
\matrix (m) [matrix of nodes, nodes in empty cells,
    column sep=4mm, row sep=3mm] { 
  \\
 1 & 2 & \textcolor{red}{4} &5  \\
   &3&  2    &  6 \\
    &   &3  &   \\
};

\draw[black, thick, -] (m-2-1) -- node[above, font=\tiny] {$1$} (m-2-2);
\draw[red, thick, -] (m-3-3) -- node[right, font=\tiny] {$\ast$}(m-2-3);
\draw[red, thick] (m-4-3) to[bend left, font=\tiny] node[left] {$\ast$} (m-2-3);
\draw  [black, thick, -] (m-2-3) --node[above, font=\tiny] {$1$}  (m-2-4);
\draw[black, thick, -] (m-3-3) --node[above, font=\tiny] {$1$}  (m-3-4);
\end{tikzpicture}\\
& \begin{tikzpicture}
\matrix (m) [matrix of nodes, nodes in empty cells,
    column sep=4mm, row sep=3mm] { 
   \\
 1 & 2 & \textcolor{red}{4} &5  \\
   &  3& 2 &   \textcolor{red} 6 \\
   &   &  &   2\\
};

\draw[black, thick, -] (m-2-1) --node[above, font=\tiny] {$1$}  (m-2-2);
\draw[red, thick, -] (m-3-3) -- node[right, font=\tiny] {$\ast$}(m-2-3);
\draw[black, thick] (m-3-2) -- node[above, font=\tiny] {$1$} (m-2-3);
\draw[black, thick, -] (m-2-3) --node[above, font=\tiny] {$1$} (m-2-4);
\draw[red, thick, -] (m-4-4) --node[right, font=\tiny] {$\ast$} (m-3-4);
\end{tikzpicture} 
\end{array}\]
\caption*{ 
For the red Set $(4,4)$ in $\mathscr{T}^{\mathcal{C}}$, the entries $2$ and $3$ both descend within $C_3(\infty)$, giving $j^{\mathcal{C}}_3 = (4,4)$ and a red $4$ appearing twice in the reverse tableau. The excluded root $x_{1,4}$ is covered by $x_{1,2}$. (section \S\ref{9.1}) For the red Set $(4,6)$ the entry $2$ is lowered into $R_2\cap C_3$ via the pair $(C_1,C_3)$, then further into $R_3\cap C_4$ via $(C_3,C_4)$. The excluded root $x_{4,6}$ is covered by $x_{4,5}$ (see \S\ref{9.2.1}).
}
\end{figure}
Example $2_4$ below illustrates the construction for the composition $(2,1,3,1,3,2,1)$, taking $C = C_3$ and $C' = C_5$.
\begin{figure}[H]
\centering
$\begin{array}{cccccccc}
 & & & C& & C'& &\\
 & C_1& C_2& C_3& C_4& C_5& C_6 & C_7 \\ \hline
R_1 & 1 & 3 & \encircled{4} & 7 & \encircled{8} & 11 & 13 \\
R_2 & 2&   & \encircled{5} &   &\encircled{9} &  12 \\
R_3 &  &  &\encircled{6}& & \encircled{10}\\
R_4
\end{array}  \xrightarrow{\,\,(C_1,C_6)\,\,} \begin{array}{cccccccc}
 & & C& & C'& &\\
  C_1& C_2& C_3& C_4& C_5& C_6 & C_7 \\ \hline
 1 & 3 & \encircled{4} & 7 & \encircled{8} & 11 & 13 \\
 2&   & \encircled{5} &   &\encircled{9} &  \red{12} \\
 \encircled{6} &  &\encircled{10}& & 12\\
\end{array} \xrightarrow{\,\,(C_4,C_7)\,\,} \begin{array}{cccccccc}
 & & C& & C'& &\\
  C_1& C_2& C_3& C_4& C_5& C_6 & C_7 \\ \hline
 1 & 3 & \encircled{4} & 7 & \encircled{8} & 11 & \red{13} \\
 2&   & \encircled{5}    &\encircled{9} &  \red{12} & 13\\
 \encircled{6} &  &\encircled{10}&  12\\
\end{array}$
\caption*{Implementing the pair $(C_1, C_6)$ of height $2$ places a red $12$ below the black $11$. Implementing $(C_4, C_8)$ of height $1$ then places a red $12$ below the black $8$, with $C'$ now of black height $1$. Yet one cannot recolor the black $8$  in red as $(C_4, C_7)$ has been used and it is not surrounded by $(C_2, C_4)$.
 }.
\end{figure}

\textbf{Example 3}.
Consider the composition $(2,1,2,1,2)$

\begin{figure}[H]
\centering
$\begin{array}{cccccc}
 & C_1& C_2& C_3& C_4& C_5\\ \hline
R_1 & 1 &3 & 4 & 6 & 7 \\
R_2 & 2 &  & 5 &   &   8
\end{array}
 \rightarrow  \begin{array}{ccccc}
 1 &3 & 4 & 6 & 7 \\
 2 &  & \textcolor{red}{5} &   &   8\\
 5
\end{array}
\begin{array}{ccc}
\nearrow& \begin{array}{ccccc} 1 &3 & \textcolor{red}{4} & 6 & 7 \\
 2 &  4& \textcolor{red}{5} &   &   8\\
 5
\end{array}& \rightarrow \begin{array}{ccccc} 1 &3 & \textcolor{red}{4} & 6 & 7 \\
 2 &  4& \textcolor{red}{5} &   &   \textcolor{red}{8}\\
 5 &    &8
\end{array}\\
\searrow & \begin{array}{ccccc} 1 &3 & 4 & \textcolor{red}{6} & 7 \\
 2 &   \textcolor{red}{5} & 6&  &   8\\
 5
\end{array}& \begin{array}{cc}
\nearrow& \begin{array}{ccccc} 1 &3 & 4 & \textcolor{red}{6} & 7 \\
 2 &   \textcolor{red}{5} & \textcolor{red}{6}&  &   8\\
 5&6
\end{array}\\
\searrow& \begin{array}{ccccc} 1 &3 & 4 & \textcolor{red}{6} & 7 \\
 2 &   \textcolor{red}{5} &{6}&  &    \textcolor{red}{8}\\
 5& &8
\end{array}
\end{array}
\end{array}$
\caption{Implement $C_1,C_3$ which recolours the black  $5$ in red.  Then implement $C_2,C_4$ which pushes this red $5$ into $R_2\cap C_2$.  Implement the pair $C_3,C_5$ is implemented. For the lowest row of the middle column, $C_3^-=C_2$ has height $2$ with its lowest entry in red $5$ and hence has black height $1$.  In this example one obtains eventually three Red Sets $(4,5,8);(5,6,6);(5,6,8)$.}
\end{figure}

 \textbf{Example 4}. Consider the composition $(1,2,3,1,1,3,2)$

\begin{figure}[H]
\centering
\[
\begin{array}{cccccccc}
 & C_1 & C_2 & C_3 & C_4 & C_5 & C_6 & C_7 \\ \hline
R_1 & 1 & 2 & 4 & 7 & 8 & 9 & 12 \\
R_2 &   & 3 & 5 &   &   & 10 & 13 \\
R_3 &   &   & 6 &   &   & 11 &
\end{array}
\quad \longrightarrow \quad\begin{array}{ccccccc}
 C_1 & C_2 & C_3 & C_4 & C_5 & C_6 & C_7 \\ \hline
 1 & 2 & 4 & \textcolor{red}{7} & \textcolor{red}{8} & 9 & 12 \\
   & 3 & 5 &   &   & 10 & \textcolor{red}{13} \\
   &   & 6 &   &   & \textcolor{red}{11} &
\end{array}\] \end{figure}
\begin{figure}[H]
\[ \quad \longrightarrow \quad
\begin{array}{ccccccc}
  C_1 & C_2 & C_3 & C_4 & C_5 & C_6 & C_7 \\ \hline
1 & 2 & 4 &{7} & {8} & 9 & 12 \\
  & 3 & 5 &   &   & 10 &\textcolor{red}{13}   \\
   & 6  & {11} &   &   & 13 &
\end{array}
\]
\caption{
  There is a component tableau $\mathscr T^\mathcal C$ with red multi-set $(7,8,11,13)$, as illustrated by the top line above.
  The contribution of the red $13$ to the excluded roots of the component tableau is illustrated in the second line.  It results from placing $13$ under $10$ and then $11$ skips over $C_4,C_5$ because they have height one, so must be skipped over to avoid a gap. This pushes $11$ under $5$ and then $6$ under $3$ creating the excluded root vector $x_{3,6}$.}
\end{figure}
\begin{figure}[H]
\centering
\[
\begin{array}{cccccccc}
 & C_1 & C_2 & C_3 & C_4 & C_5 & C_6 & C_7 \\ \hline
R_1 & 1 & 2 & 4 & 7 & 8 & 9 & 12 \\
R_2 &   & 3 & 5 &   &   & 10 & 13 \\
R_3 &   &   & 6 &   &   & 11 &
\end{array}
\quad \longrightarrow \quad\begin{array}{ccccccc}
 C_1 & C_2 & C_3 & C_4 & C_5 & C_6 & C_7 \\ \hline
 1 & 2 & 4 & \textcolor{red}{7} & {8} & 9 & 12 \\
 3  & 5 & 7 &   &   & 10 & {13} \\
   &  6 &  &   &   &{11} &
\end{array}\] \end{figure}
\begin{figure}[H]
\[ \quad \longrightarrow \quad
\begin{array}{ccccccc}
  C_1 & C_2 & C_3 & C_4 & C_5 & C_6 & C_7 \\ \hline
 1 & 2 & 4 & \textcolor{red}{7} & \textcolor{red}{8} & 9 & 12  \\
  3& 5 & 7 & 8  &   & 10 & {13}   \\
   & 6  & &   &   & 11 &
\end{array}\quad \longrightarrow \quad
\begin{array}{ccccccc}
  C_1 & C_2 & C_3 & C_4 & C_5 & C_6 & C_7 \\ \hline
 1 & 2 & 4 & \textcolor{red}{7} & \textcolor{red}{8} & 9 & 12  \\
  3& 5 & 7 & 8  &   & 10 & {13}   \\
   & 6  & &   &   &  \textcolor{red}{11} &   \\
      & 11  & &   &   &  &
\end{array}
\] \[\quad \longrightarrow \quad
\begin{array}{ccccccc}
  C_1 & C_2 & C_3 & C_4 & C_5 & C_6 & C_7 \\ \hline
 1 & 2 & 4 & \textcolor{red}{7} & \textcolor{red}{8} & 9 & 12  \\
  3& 5 & 7 & 8  &   & 10 & \textcolor{red}{13}   \\
   6&  \textcolor{red}{11} & &   &   &  {13} &   \\
     11 &   & &   &   &  &
\end{array}\]
\caption*{
 Now consider the partial reverse tableau corresponding to the subset $(7,8,11)$ of the red multi-set $(7,8,11,13)$ above.
 In this  $C_3$ (resp. $C_4$) have height $2$ with entries $4,7$ (resp. red $7,8$), whilst $C_2$ has height $4$ with entries $2,5,6,11$.  In this $C_1$ has height $2$ obtaining its lowest entry $3$ by translation.  Redressing by our slightly modified rule replaces $C_2^{>2}$ by $C_6^{>2}$ and places the former under $C_1$, so in particular placing $6$ under $3$. Finally to obtain the reverse tableau with red multi-set $(7,8,11,13)$, the black $13$ is placed in $R_3\cap C_5$ under $10$ recolouring the former in red and the red $11$ is placed in $R_3\cap C_2$ under $5$.  Without redressing this is equivalent to the red $11$ skipping over $C_3,C_4$ even though they have height $2$ and pushing $C_2^{>2}$ in $C_1$ under $3$.  This is the extreme shifting referred to in \ref {4.3.2}. Thus the reverse tableau admits  the excluded root vector $x_{3,6}$, \textit{which would have been absent without extreme shifting}.
 }
\end{figure}


\textbf{Example 5}.
Consider the composition $(3,4,2,1,2,4,3,1)$

\begin{figure}[H]
\centering
\[
\begin{array}{ccccccccc}
 & C_1 & C_2 & C_3 & C_4 & C_5 & C_6 & C_7 & C_8 \\ \hline
R_1 & 1 & 4 & 8 & 10 & 11 & 13 & 17 & 20 \\
R_2 & 2 & 5 & 9 &   & 12 & 14 & 18 &   \\
R_3 &  3 &6 &   &   &   &    15  &   19\\
R_4 &   &  7 &   &   &   &     16 &
\end{array} \quad \longrightarrow \quad \begin{array}{ cccccccc}
  C_1 & C_2 & C_3 & C_4 & C_5 & C_6 & C_7 & C_8 \\ \hline
  1 & 4 & 8 & 10 & {11} & 13 & 17 & 20 \\
  2 & 5 & 9 &   & \textcolor{red}{12} & 14 & 18 &   \\
   3 &6 & 12  &   &   &   {15}  &   19&\\
    &  7 &   &   &   &    {16} &
\end{array}
\]
\[\quad \longrightarrow \quad \begin{array}{ cccccccc}
  C_1 & C_2 & C_3 & C_4 & C_5 & C_6 & C_7 & C_8 \\ \hline
  1 & 4 & 8 & 10 &  \textcolor{red}{11} & 13 & 17 & 20 \\
  2 & 5 & 9 &  11 &  \textcolor{red}{12} & 14 & 18 &   \\
   3 &6 & 12  &   &   &   {15}  &   19&\\
    &  7 &   &   &   &     {16} &   \\

\end{array}\quad \longrightarrow \quad \begin{array}{ cccccccc}
  C_1 & C_2 & C_3 & C_4 & C_5 & C_6 & C_7 & C_8 \\ \hline
  1 & 4 & 8 & 10 & \textcolor{red}{11} & 13 & 17 & 20 \\
  2 & 5 & 9 &  11 &\textcolor{red} {12} & 14 & 18 &   \\
   3 &6 & 12  &   &   &    {15}  &   19&\\
    &  7 &   &   &   &     \textcolor{red}{16} &   \\
    &16
\end{array} \]\end{figure}
\begin{figure}[H] \[\quad \longrightarrow \quad \begin{array}{ cccccccc}
  C_1 & C_2 & C_3 & C_4 & C_5 & C_6 & C_7 & C_8 \\ \hline
  1 & 4 & 8 & 10 & \textcolor{red}{11} & 13 & 17 & 20 \\
  2 & 5 & 9 & 11  & \textcolor{red}{12} & 14 & 18 &   \\
   3 &6 &  12 &   &   &     \textcolor{red}{15}  &   19&\\
    & 7 & 15  &   &   &     \textcolor{red}{16} &   \\
    &16&
\end{array}\]
\caption{Let the $10$-string go down by two consecutive steps of height two. This gives the Red Set $(11,12,15,16)$ as shown in the first line.
 To construct the reverse tableau, recall that according to our ordering of the red multi-set we implement in the order $12,11,16,15$.
  First $12$ goes under $9$ and then $11$ under $10$.  Then $16$ goes under $7$ finally $15$ under $12$.
  Observe that the amalgamated column $\hat{C_5}$ of \cite [4.3.6]{FJ5} which has entries $10,11,12,15,16$ appears in the reverse tableau skewed to the left.}
\end{figure}
An alternative reverse tableau with the red set $11,12,15,16$ obtained by the shifting $R^4 \cap C_2 C_3$ to the left. The matrix following it shows that the additional excluded root are covered in the sense of section  \S\ref{9.2}
\[
 \begin{array}{ cccccccc}
  C_1 & C_2 & C_3 & C_4 & C_5 & C_6 & C_7 & C_8 \\ \hline
  1 & 4 & 8 & 10 & \textcolor{red}{11} & 13 & 17 & 20 \\
  2 & 5 & 9 &  11 & \textcolor{red}{12} & 14 & 18 &   \\
   3 &6 &  12 &   &   &    \textcolor{red}{15}  &   19&\\
    7&  15 &   &   &   &     \textcolor{red}{16} &   \\
    16
\end{array}
\]
\begin{figure}[H]
\begin{tikzpicture}\tiny
\matrix[matrix of math nodes, left delimiter=(, right delimiter=)] (m) {
1 & 0 & 0 & x_{1,4} & x_{1,5} & x_{1,6} & \textcolor{blue}{x_{1,7}} & x_{1,8} & x_{1,9} & x_{1,10} & x_{1,11} & x_{1,12} & x_{1,13} & x_{1,14} & x_{1,15} &\textcolor{blue}{x_{1,16}} & x_{1,17} & x_{1,18} & x_{1,19} & x_{1,20} \\
0 & 1 & 0 & x_{2,4} & x_{2,5} & x_{2,6} & \textcolor{blue}{x_{2,7}} & x_{2,8} & x_{2,9} & x_{2,10} & x_{2,11} & x_{2,12} & x_{2,13} & x_{2,14} & x_{2,15} & \textcolor{blue}{x_{2,16} }& x_{2,17} & x_{2,18} & x_{2,19} & x_{2,20} \\
0 & 0 & 1 & x_{3,4} & x_{3,5} & x_{3,6} & \textcolor{blue}{x_{3,7}} & x_{3,8} & x_{3,9} & x_{3,10} & x_{3,11} & x_{3,12} & x_{3,13} & x_{3,14} & x_{3,15} & \textcolor{blue}{x_{3,16} }& x_{3,17} & x_{3,18} & x_{3,19} & x_{3,20} \\
0 & 0 & 0 & 1 & 0 & 0 & 0 & 0 & x_{4,9} & x_{4,10} & x_{4,11} & x_{4,12} & x_{4,13} & x_{4,14} &\textcolor{blue}{x_{4,15} }& \textcolor{blue}{x_{4,16} }& x_{4,17} & x_{4,18} & x_{4,19} & x_{4,20} \\
0 & 0 & 0 & 0 & 1 & 0 & 0 & 0 & x_{5,9} & x_{5,10} & x_{5,11} & x_{5,12} & x_{5,13} & x_{5,14} & \textcolor{blue}{x_{5,15} }&\textcolor{blue}{x_{5,16}}& x_{5,17} & x_{5,18} & x_{5,19} & x_{5,20} \\
0 & 0 & 0 & 0 & 0 & 1 & 0 & 0 & x_{6,9} & x_{6,10} & x_{6,11} & x_{6,12} & x_{6,13} & x_{6,14} & \textcolor{blue}{x_{6,15}} &\textcolor{blue}{x_{6,16}} & x_{6,17} & x_{6,18} & x_{6,19} & x_{6,20} \\
0 & 0 & 0 & 0 & 0 & 0 & 1 & 0 & x_{7,9} & x_{7,10} & x_{7,11} & x_{7,12} & x_{7,13} & x_{7,14} & x_{7,15} & \textcolor{blue}{x_{7,16}} & x_{7,17} & x_{7,18} & x_{7,19} & x_{7,20} \\
0 & 0 & 0 & 0 & 0 & 0 & 0 & 1 & 0 & x_{8,10} & x_{8,11} & \textcolor{blue}{x_{8,12} }& x_{8,13} & x_{8,14} &\textcolor{blue}{x_{8,15}} & \textcolor{blue}{x_{8,16}} & x_{8,17} & x_{8,18} & x_{8,19} & x_{8,20} \\
0 & 0 & 0 & 0 & 0 & 0 & 0 & 0 & 1 & x_{9,10} & x_{9,11}& \textcolor{blue}{x_{9,12}} & x_{9,13} & x_{9,14} &\textcolor{blue}{x_{9,15} }& \textcolor{blue}{x_{9,16}} & x_{9,17} & x_{9,18} & x_{9,19} & x_{9,20} \\
0 & 0 & 0 & 0 & 0 & 0 & 0 & 0 & 0 & 1 & \textcolor{red}{\ast} & \textcolor{red}{\ast} & x_{10,13} & x_{10,14} &  \textcolor{red}{\ast} &  \textcolor{red}{\ast}& x_{10,17} & x_{10,18} & x_{10,19} & x_{10,20} \\
0 & 0 & 0 & 0 & 0 & 0 & 0 & 0 & 0 & 0 & 1 & 0 & x_{11,13} & x_{11,14} &\textcolor{blue}{x_{11,15} }& \textcolor{blue}{x_{11,16}} & x_{11,17} & x_{11,18} & x_{11,19} & x_{11,20} \\
0 & 0 & 0 & 0 & 0 & 0 & 0 & 0 & 0 & 0 & 0 & 1 & x_{12,13} & x_{12,14} &\textcolor{blue}{x_{12,15}} & \textcolor{blue}{x_{12,16}} & x_{12,17} & x_{12,18} & x_{12,19} & x_{12,20} \\
0 & 0 & 0 & 0 & 0 & 0 & 0 & 0 & 0 & 0 & 0 & 0 & 1 & 0 & 0 & 0 & x_{13,17} & x_{13,18} & x_{13,19} & x_{13,20} \\
0 & 0 & 0 & 0 & 0 & 0 & 0 & 0 & 0 & 0 & 0 & 0 & 0 & 1 & 0 & 0 & x_{14,17} & x_{14,18} & x_{14,19} & x_{14,20} \\
0 & 0 & 0 & 0 & 0 & 0 & 0 & 0 & 0 & 0 & 0 & 0 & 0 & 0 & 1 & 0 & x_{15,17} & x_{15,18} & x_{15,19} & x_{15,20} \\
0 & 0 & 0 & 0 & 0 & 0 & 0 & 0 & 0 & 0 & 0 & 0 & 0 & 0 & 0 & 1 & x_{16,17} & x_{16,18} & x_{16,19} & x_{16,20} \\
0 & 0 & 0 & 0 & 0 & 0 & 0 & 0 & 0 & 0 & 0 & 0 & 0 & 0 & 0 & 0 & 1 & 0 & 0 & x_{17,20} \\
0 & 0 & 0 & 0 & 0 & 0 & 0 & 0 & 0 & 0 & 0 & 0 & 0 & 0 & 0 & 0 & 0 & 1 & 0 & x_{18,20} \\
0 & 0 & 0 & 0 & 0 & 0 & 0 & 0 & 0 & 0 & 0 & 0 & 0 & 0 & 0 & 0 & 0 & 0 & 1 & x_{19,20} \\
0 & 0 & 0 & 0 & 0 & 0 & 0 & 0 & 0 & 0 & 0 & 0 & 0 & 0 & 0 & 0 & 0 & 0 & 0 & 1 \\
};

\draw[thick, red] (m-1-1.north west) rectangle (m-3-3.south east);   
\draw[thick, blue] (m-4-4.north west) rectangle (m-7-7.south east);  
\draw[thick, green] (m-8-8.north west) rectangle (m-9-9.south east); 
\draw[thick, orange] (m-10-10.north west) rectangle (m-10-10.south east); 
\draw[thick, green] (m-11-11.north west) rectangle (m-12-12.south east); 
\draw[thick, blue] (m-13-13.north west) rectangle (m-16-16.south east); 
\draw[thick,red] (m-17-17.north west) rectangle (m-19-19.south east); 
\draw[thick, orange] (m-20-20.north west) rectangle (m-20-20.south east); 
\end{tikzpicture}
\caption*{The elements shown in 
{blue} represent the excluded roots.
 the excluded roots ofor a maximal pushing to the left.
Nonetheless, the action of $P$ on the non-excluded roots of the component and reverse tableaux determines the same geometric component, which in this case corresponds to an orbital variety.
 }

\end{figure}

\textbf{Example 6}.
Consider the composition $(2,1,1,1,2)$

\begin{figure}[H]
\centering
$\begin{array}{cccccc}
 & C_1& C_2& C_3& C_4& C_5\\ \hline
R_1 & 1 & 3 & {4} & 5 & 6 \\
R_2 & 2 &  &  &   & 7  \\
\end{array}\longrightarrow \begin{array}{ccccc}
 C_1& C_2& C_3& C_4& C_5\\ \hline
 1 & 3 & \textcolor{red}{4} & \textcolor{red}{5} & 6 \\
 2 & 4 &  5&   & \textcolor{red}{7}  \\
 & &7
\end{array}$

\begin{tikzpicture}
\matrix[matrix of math nodes, left delimiter=(, right delimiter=)] (m) {
1 & 0 & x_{1,3} & x_{1,4} & x_{1,5} & x_{1,6} & x_{1,7} \\
0 & 1 & x_{2,3} & x_{2,4} & x_{2,5} & x_{2,6} & x_{2,7} \\
0 & 0 & 1 & \textcolor{red}{\ast} & x_{3,5} & x_{3,6} & x_{3,7} \\
0 & 0 & 0 & 1 &\textcolor{red}{\ast} & x_{4,6} & \textcolor{red}{\ast} \\
0 & 0 & 0 & 0 & 1 & x_{5,6} & \textcolor{blue}{x_{5,7}} \\
0 & 0 & 0 & 0 & 0 & 1 & 0 \\
0 & 0 & 0 & 0 & 0 & 0 & 1 \\
};

\draw[thick, purple]    (m-1-1.north west) rectangle (m-2-2.south east);  
\draw[thick, blue]   (m-3-3.north west) rectangle (m-3-3.south east);  
\draw[thick, green]  (m-4-4.north west) rectangle (m-4-4.south east);  
\draw[thick, orange] (m-5-5.north west) rectangle (m-5-5.south east);  
\draw[thick, purple] (m-6-6.north west) rectangle (m-7-7.south east);  

\end{tikzpicture}

\caption{A reverse tableau with red multi-set ${4,5,7}$ is constructed as follows. First a black $4$ is placed under $3$ and the original black $4$ is recoloured in red. Then a black $5$ is placed under this red $4$ and the oribinal black $5$ recoloured in red. Finally a black $7$ is placed under this black $5$ and the original balck $7$ recoloured.  This was all done simultaneously in the first lime of figure. It produces an excluded root indicated in \textbf{M} in blue. In the component tableau the red line $\ell_{4,7}$ also produces this excluded root, which is covered by the line $\ell_{5,6}$ labelled by $1$.}
\end{figure}

\begin{figure}[H]
\centering
$\begin{array}{cccccc}
 & C_1& C_2& C_3& C_4& C_5\\ \hline
R_1 & 1 & 3 & {4} & 5 & 6 \\
R_2 & 2 &  &  &   & 7  \\
\end{array}\longrightarrow \begin{array}{ccc}
\nearrow & \begin{array}{ccccc}
 C_1& C_2& C_3& C_4& C_5\\ \hline
 1 & 3 & \textcolor{red}{4} & \textcolor{red}{5} & 6 \\
 2 & 4 &   \textcolor{red}{5}&   &{7}  \\
  &5 &
\end{array}&\longrightarrow
\begin{tikzpicture}
\matrix[matrix of math nodes, left delimiter=(, right delimiter=)] (m) {
1 & 0 & x_{1,3} & x_{1,4} & {x_{1,5}} & x_{1,6} & x_{1,7} \\
0 & 1 & x_{2,3} & x_{2,4} & {x_{2,5}} & x_{2,6} & x_{2,7} \\
0 & 0 & 1 & \textcolor{red}{\ast} & \textcolor{red}{\ast} & x_{3,6} & x_{3,7} \\
0 & 0 & 0 & 1 &\textcolor{red}{\ast} & x_{4,6} & x_{4,7} \\
0 & 0 & 0 & 0 & 1 & x_{5,6} & {x_{5,7}} \\
0 & 0 & 0 & 0 & 0 & 1 & 0 \\
0 & 0 & 0 & 0 & 0 & 0 & 1 \\
};
\draw[thick, purple]    (m-1-1.north west) rectangle (m-2-2.south east);  
\draw[thick, blue]   (m-3-3.north west) rectangle (m-3-3.south east);  
\draw[thick, green]  (m-4-4.north west) rectangle (m-4-4.south east);  
\draw[thick, orange] (m-5-5.north west) rectangle (m-5-5.south east);  
\draw[thick, purple] (m-6-6.north west) rectangle (m-7-7.south east);  
\end{tikzpicture}\\
\searrow & \begin{array}{ccccc}
 C_1& C_2& C_3& C_4& C_5\\ \hline
 1 & 3 & \textcolor{red}{4} & \textcolor{red}{5} & 6 \\
 2 & 4 &   \textcolor{red}{5}&   &{7}  \\
 5 & &
\end{array}&\longrightarrow
\begin{tikzpicture}
\matrix[matrix of math nodes, left delimiter=(, right delimiter=)] (m) {
1 & 0 & x_{1,3} & x_{1,4} & \textcolor{blue}{x_{1,5}} & x_{1,6} & x_{1,7} \\
0 & 1 & x_{2,3} & x_{2,4} & \textcolor{blue}{x_{2,5}} & x_{2,6} & x_{2,7} \\
0 & 0 & 1 & \textcolor{red}{\ast} & \textcolor{red}{\ast}& x_{3,6} & x_{3,7} \\
0 & 0 & 0 & 1 &\textcolor{red}{\ast} & x_{4,6} & x_{4,7} \\
0 & 0 & 0 & 0 & 1 & x_{5,6} & {x_{5,7}} \\
0 & 0 & 0 & 0 & 0 & 1 & 0 \\
0 & 0 & 0 & 0 & 0 & 0 & 1 \\
};
\draw[thick, purple]    (m-1-1.north west) rectangle (m-2-2.south east);  
\draw[thick, blue]   (m-3-3.north west) rectangle (m-3-3.south east);  
\draw[thick, green]  (m-4-4.north west) rectangle (m-4-4.south east);  
\draw[thick, orange] (m-5-5.north west) rectangle (m-5-5.south east);  
\draw[thick, purple] (m-6-6.north west) rectangle (m-7-7.south east);  
\end{tikzpicture}

\end{array}$
\caption*{Take the composition $(2,1,1,1,2)$. Consider the Red Set $4,5,5$.  The set $X(\ell^\mathcal C_{4,5})$ is defined by placing $5$ under $4$ in the component tableau.  This is also true in the reverse tableau except that $4$ is now in red.  Thus in computing $X(\ell_{4,5}^{\psi(\mathcal C)})$ we must allow red to be a colour of the first entry.}
\end{figure}

\textbf{Example 7}.
Consider the composition $(1,2,1,2,1)$.

\begin{figure}[H]
\centering
$\begin{array}{cccccc}
 & C_1& C_2& C_3& C_4& C_5\\ \hline
R_1 & 1 & 2 & 4 & 5 & 7 \\
R_2 &   & 3 &  & 6 &   \\
\end{array} \longrightarrow $ $\begin{array}{ccccc}
  C_1& C_2& C_3& C_4& C_5\\ \hline
 1 & 2 &  \textcolor{red}{4} & 5 & \textcolor{red}{7} \\
  3 &  \textcolor{red}{4} & 6 & 7 &   \\
  4
\end{array} $ $\longrightarrow \begin{tikzpicture}
\matrix[matrix of math nodes, left delimiter=(, right delimiter=)] (m) {
1 & x_{1,2} & \textcolor{blue}{x_{1,3} }& \textcolor{blue}{x_{1,4} }& x_{1,5} & x_{1,6} & x_{1,7} \\
0 & 1 & 0 & \textcolor{red}{\ast} & x_{2,5} & x_{2,6} & x_{2,7} \\
0 & 0 & 1 &  \textcolor{red}{\ast}& x_{3,5} & x_{3,6} & x_{3,7} \\
0 & 0 & 0 & 1 & x_{4,5} & \textcolor{blue}{x_{4,6} }& x_{4,7} \\
0 & 0 & 0 & 0 & 1 & 0 & \textcolor{red}{\ast} \\
0 & 0 & 0 & 0 & 0 & 1 & x_{6,7} \\
0 & 0 & 0 & 0 & 0 & 0 & 1 \\
};

\draw[thick, red]    (m-1-1.north west) rectangle (m-1-1.south east);   
\draw[thick, orange]   (m-2-2.north west) rectangle (m-3-3.south east);   
\draw[thick, green]  (m-4-4.north west) rectangle (m-4-4.south east);   
\draw[thick, orange] (m-5-5.north west) rectangle (m-6-6.south east);   
\draw[thick, purple] (m-7-7.north west) rectangle (m-7-7.south east);   
\end{tikzpicture}
$
\caption{This has a component with red multi-set $(4,4,7)$. In the corresponding reverse tableau, a (black) $6$ lies strictly to the left of the black $4$.  On the other hand this black $6$ lies strictly below a red $4$ in the same column,  so only by virtue of the latter is $x_{4,6}$ an excluded root. Again in this example $x_{1,4}$ lies in $X(\mathscr R^{\psi(\mathcal C)})$ but not in $X(\mathscr T^\mathcal C)$.
}
\end{figure}

\textbf{Example 8}.
Consider the composition $(2,1,1,2,2)$.

\begin{figure}[H]
\centering
$\begin{array}{cccccc}
 & C_1& C_2& C_3& C_4& C_5\\ \hline
R_1 & 1 & 3 & 4 & 5 & 7 \\
R_2 & 2 &  &   &   6&  8 \\

\end{array} \longrightarrow \begin{array}{cccccc}
 C_1& C_2& C_3& C_4& C_5\\ \hline
 1 & 3 & \textcolor{red}{4} & 5 & 7 \\
 2 & 4 &   &   \textcolor{red}{6}&  \textcolor{red}{8} \\
   &  6 &   & 8  &   \\
\end{array}$\\

$\begin{tikzpicture}[
]
\matrix (m) [matrix of nodes, nodes in empty cells,
    column sep=0.5mm, row sep=0.5mm] {
$C_1$ & $C_2$ & $C_3$ & $C_4$ & $C_5$ \\\hline 
1 & 3 & \textcolor{red}{4} & 5 & 7 \\
2 & 4 &  & \textcolor{red}{6} & \textcolor{red}{8} \\
  & 6 &  & 8 &  \\
};

\draw[blue, thick, -] (m-2-2) -- (m-4-4);

\end{tikzpicture}$
\caption{The component with red multi-set $4,6,8$.  Here  $1$ appears in row $3$ in column $8$ of \textbf{M} following  copies of $\ast$ in columns $4,6$. By the requirement of \ref{8.3}, in the reverse tableau the right going line $\ell_{3,8}$ labelled by a $1$ must go from the singleton entry of $3$ to the black entry of $8$ and thern down by two rows. This requiremnt was needed for \ref {8.3}${(*)}$ which is indeed satisfied in this example.  The remaining four lines in the reverse tableau labelled by a $1$ are horizontal.}
\end{figure}

\textbf{Example 9}.
Consider the composition $(3,2,1,3,2,1,2)$.

\begin{figure}[H]
\centering
$\begin{array}{cccccccc}
 & C_1& C_2& C_3& C_4& C_5& C_6 & C_7 \\ \hline
R_1 & 1 & 4 & 6 & 7 & 10 & 12 & 13 \\
R_2 & 2& 5& & 8 &  11 & &  14 \\
R_3 & 3& & &9\\
R_4
\end{array} \longrightarrow \begin{array}{cccccccc}
  C_1& C_2& C_3& C_4& C_5& C_6 & C_7 \\ \hline
  1 & 4 & 6 & 7 & 10 & \textcolor{red}{12} & 13 \\
  2& 5&  \textcolor{red}{8} &  11 &\textcolor{red}{12} &  &14 \\
 3& 8&  \textcolor{red}{9}&{12} \\
 9& & \\
\end{array}$

\caption{The component with red multi-set $8,9,12,12$.  The line $\ell_{9,12}$ goes up by $m=2$ rows in the component tableau whilst the new black $12$ goes down by $m=2$ rows in the reverse tableau.  Thus the line $\ell_{9,12}$ from the rightmost $9$ to the leftmost $12$, is horizontal in the reverse tableau.}
\end{figure}

\textbf{Example $10$}. This illustrates the factorisation explained in \ref {9.5} to recover all the reverse tableaux corresponding to the red multi-sets of Example $1$.
%

\begin{figure}[H]
 \[\begin{array}{ ccccccccc}
 &C_1 & C_2 & C_3 & C_4 & C_5 & C_6 & C_7 \\ \hline
  R_1 &1&2&4&7&8&9&12\\
    R_2& &3&5& &&10&13\\
    R_3&& & 6& &&11&
  \end{array} \quad \longrightarrow \quad  \begin{array}{ ccccccccc}
 C_1 & C_2 & C_3 & C_4 & C_5 & C_6 & C_7 \\ \hline
1&2&4&\textcolor{red}{7}&8&9&12\\
3 &5&7& &&10&13\\
& 6& & &&11&
  \end{array}\quad \longrightarrow \quad \]
  \end{figure}

  \begin{figure}[H]
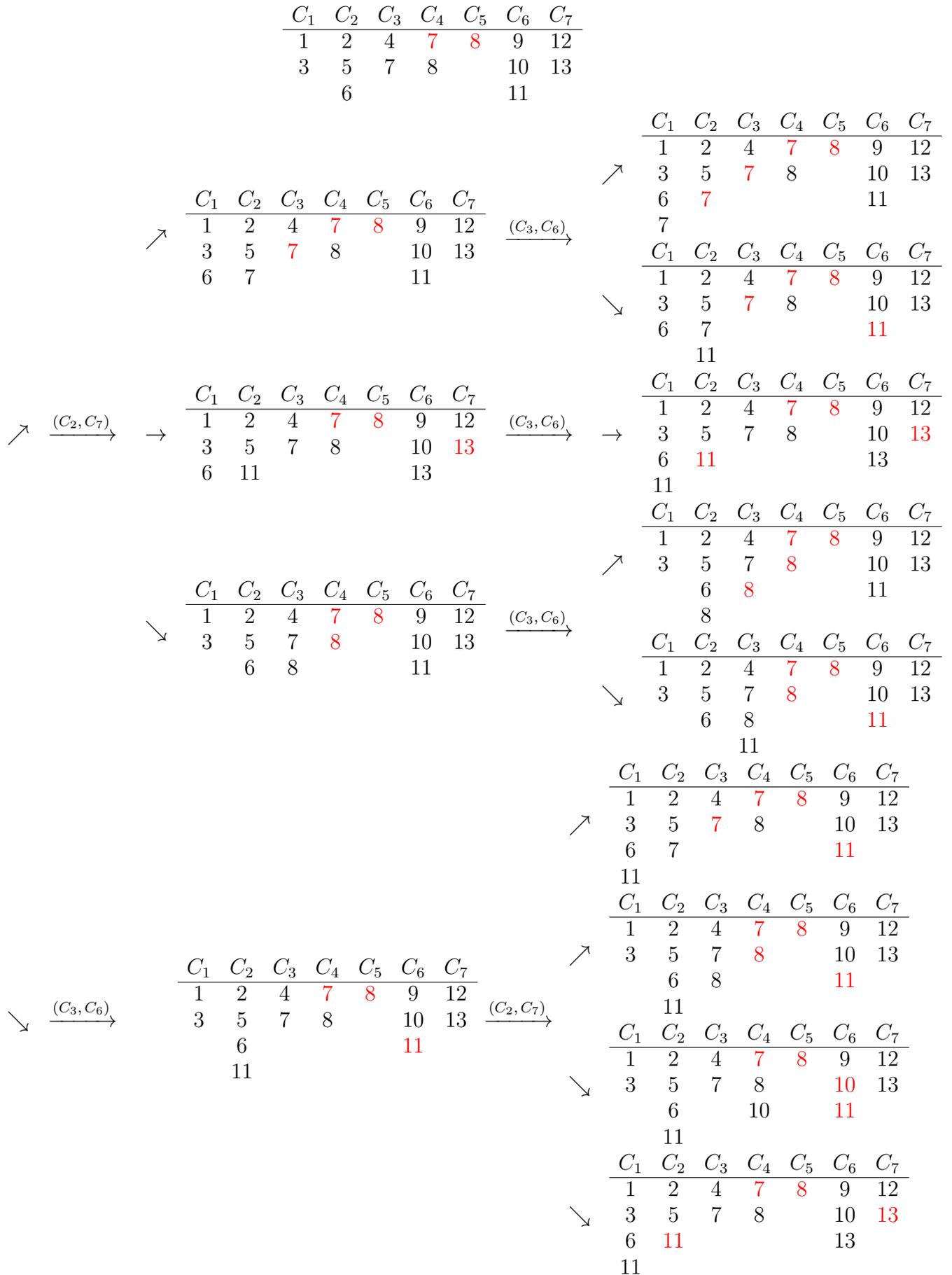

  $ \begin{array}{ ccccccccc}
 C_1 & C_2 & C_3 & C_4 & C_5 & C_6 & C_7 \\ \hline
1&2&4&\textcolor{red}{7}&\textcolor{red}{8}&9&12\\
3 &5&7&8 &&10&13\\
& 6& & &&11&
  \end{array} $
$\begin{array}{ ccc} \nearrow&\xrightarrow{\text{$(C_2,C_7)$}} &\begin{array}{cccc} \nearrow&\begin{array}{ ccccccccc}
 C_1 & C_2 & C_3 & C_4 & C_5 & C_6 & C_7 \\ \hline
1&2&4&\textcolor{red}{7}&\textcolor{red}{8}&9&12\\
3 &5&\textcolor{red}{7}&8 &&10&13\\
6& 7& & &&11&
  \end{array} &\xrightarrow{\text{$(C_3,C_6)$}} &\begin{array}{cc}  \nearrow&\begin{array}{ ccccccccc}
 C_1 & C_2 & C_3 & C_4 & C_5 & C_6 & C_7 \\ \hline
1&2&4&\textcolor{red}{7}&\textcolor{red}{8}&9&12\\
3 &5&\textcolor{red}{7}&8 &&10&13\\
6&\textcolor{red}{7}& & &&11&\\
{7}
  \end{array}\\
  \searrow&\begin{array}{ ccccccccc}
 C_1 & C_2 & C_3 & C_4 & C_5 & C_6 & C_7 \\ \hline
1&2&4&\textcolor{red}{7}&\textcolor{red}{8}&9&12\\
3 &5&\textcolor{red}{7}&8 &&10&13\\
6&{7}& & &&\textcolor{red}{11}&\\
&11
  \end{array}
    \end{array}\\
   \rightarrow&\begin{array}{ ccccccccc}
 C_1 & C_2 & C_3 & C_4 & C_5 & C_6 & C_7 \\ \hline
1&2&4&\textcolor{red}{7}&\textcolor{red}{8}&9&12\\
3 &5&{7}&8 &&10&\textcolor{red}{13}\\
6 & 11&&& &13&
  \end{array} &\xrightarrow{\text{$(C_3,C_6)$}} &\begin{array}{cc}  \rightarrow&\begin{array}{ ccccccccc}
 C_1 & C_2 & C_3 & C_4 & C_5 & C_6 & C_7 \\ \hline
1&2&4&\textcolor{red}{7}&\textcolor{red}{8}&9&12\\
3 &5&{7}&8 &&10&\textcolor{red}{13}\\
6&\textcolor{red}{11}& & &&{13}&\\
11
  \end{array}\end{array}\\
   \searrow&\begin{array}{ ccccccccc}
 C_1 & C_2 & C_3 & C_4 & C_5 & C_6 & C_7 \\ \hline
1&2&4&\textcolor{red}{7}&\textcolor{red}{8}&9&12\\
3 &5&{7}&\textcolor{red}{8} &&10&13\\
&6 & 8& &&11&
  \end{array}& \xrightarrow{\text{$(C_3,C_6)$}} &\begin{array}{cc}  \nearrow&\begin{array}{ ccccccccc}
 C_1 & C_2 & C_3 & C_4 & C_5 & C_6 & C_7 \\ \hline
1&2&4&\textcolor{red}{7}&\textcolor{red}{8}&9&12\\
3 &5&{7}&\textcolor{red}{8} &&10&13\\
&6 & \textcolor{red}{8} & &&11&\\
&8
  \end{array}\\
  \searrow&\begin{array}{ ccccccccc}
 C_1 & C_2 & C_3 & C_4 & C_5 & C_6 & C_7 \\ \hline
1&2&4&\textcolor{red}{7}&\textcolor{red}{8}&9&12\\
3 &5&{7}&\textcolor{red}{8} &&10&13\\
&6 &{8} & && \textcolor{red}{11}&\\
& &11
  \end{array}
  \end{array}\\

  \end{array} \\
 \searrow&\xrightarrow{\text{$(C_3,C_6)$}}& \begin{array}{ ccccccccc}
 C_1 & C_2 & C_3 & C_4 & C_5 & C_6 & C_7 \\ \hline
1&2&4&\textcolor{red}{7}&\textcolor{red}{8}&9&12\\
3 &5&{7}&8 &&10&13\\
& 6& & &&\textcolor{red} {11}&\\
&11
  \end{array} \xrightarrow{\text{$(C_2,C_7)$}} \begin{array}{cc}  \nearrow&\begin{array}{ ccccccccc}
 C_1 & C_2 & C_3 & C_4 & C_5 & C_6 & C_7 \\ \hline
1&2&4&\textcolor{red}{7}&\textcolor{red}{8}&9&12\\
3 &5&\textcolor{red}{7}&8&&10&13\\
6& 7& & &&\textcolor{red} {11}&\\
11&
  \end{array}\\
   \nearrow&\begin{array}{ ccccccccc}
 C_1 & C_2 & C_3 & C_4 & C_5 & C_6 & C_7 \\ \hline
1&2&4&\textcolor{red}{7}&\textcolor{red}{8}&9&12\\
3 &5&{7}&\textcolor{red}{8} &&10&13\\
& 6&8 & &&\textcolor{red} {11}&\\
&11
  \end{array}\\
   \searrow&\begin{array}{ ccccccccc}
 C_1 & C_2 & C_3 & C_4 & C_5 & C_6 & C_7 \\ \hline
1&2&4&\textcolor{red}{7}&\textcolor{red}{8}&9&12\\
3 &5&{7}&8 &&\textcolor{red} {10}&13\\
& 6& &10 &&\textcolor{red} {11}&\\
&11
  \end{array}\\
   \searrow&\begin{array}{ ccccccccc}
 C_1 & C_2 & C_3 & C_4 & C_5 & C_6 & C_7 \\ \hline
1&2&4&\textcolor{red}{7}&\textcolor{red}{8}&9&12\\
3 &5&{7}&8 &&10&\textcolor{red} {13}\\
 6& \textcolor{red} {11}& &&&13&\\
11
  \end{array}
  \end{array}
      \end{array}
$
\caption{The flow chart for the two factorisations given in  \ref {9.5}}.

  \end{figure}

\end{document}